\documentclass[twocolumn]{autart}

\usepackage{amsmath} 
\usepackage{amssymb,latexsym}
\usepackage{graphicx,color}
\usepackage{floatflt}

\usepackage{epstopdf}
\usepackage{float}
\usepackage{subcaption}
\usepackage{flushend}
\usepackage{mathrsfs}
\usepackage{cases}
\usepackage[overload]{empheq}
\usepackage[retainorgcmds]{IEEEtrantools}
\usepackage{enumitem}

\renewcommand\tilde{\widetilde}

\newtheorem{nnassumption}{\bf Assumption}
\newenvironment{assumption}{\begin{nnassumption}\it}{\end{nnassumption}}
\newtheorem{nntheorem}{\bf Theorem}
\newenvironment{theorem}{\begin{nntheorem}\it}{\end{nntheorem}}
\newtheorem{nncorollary}{\bf Corollary}

\newtheorem{nndefinition}{\bf Definition}
\newenvironment{definition}{\begin{nndefinition}\it}{\end{nndefinition}}
\newtheorem{nnproposition}{\bf Proposition}
\newenvironment{proposition}{\begin{nnproposition}\it}{\end{nnproposition}}
\newtheorem{nnproblem}{\bf Problem}

\newtheorem{nnlemma}{\bf Lemma}
\newenvironment{lemma}{\begin{nnlemma}\it}{\end{nnlemma}}
\newtheorem{nnremark}{\bf Remark}
\newenvironment{remark}{\begin{nnremark} \rm }{\hfill \hspace*{1pt}\hfill $\Box$\end{nnremark}}
\newtheorem{nnexample}{\bf Example}

\newcommand{\vv{\vspace{2em}}}

\newcommand\RR{\mathbb{R}}
\begin{document}

\begin{frontmatter}

\title{ On-line detection of qualitative dynamical changes in nonlinear systems: the resting-oscillation case\thanksref{footnoteinfo}}

\thanks[footnoteinfo]{A preliminary version of this work was presented at the 20th IFAC World Congress, Toulouse, France, 2017 \cite{tangetalifac17w}. This work was supported by the ANR under the grant SEPICOT (ANR 12 JS03 004 01), by the French ``R{\'e}gion Grand Est" through a fellowship grant 2016-2017, and by DGAPA-UNAM under the grant PAPIIT RA105518.}
\author[cristal]{Ying Tang} 
\author[unam]{Alessio Franci} 
\author[cran]{Romain Postoyan}

\address[cristal]{CRIStAL CNRS UMR 9189, Universit{\'e} de Lille, France {\tt ying.tang@univ-lille.fr}}
\address[unam]{Department of Mathematics, Universidad National Aut{\'o}noma de M{\'e}xico, Mexico {\tt afranci@ciencias.unam.mx}}
\address[cran]{Universit{\'e} de Lorraine, CNRS, CRAN, F-54000 Nancy, France {\tt romain.postoyan@univ-lorraine.fr}}

\maketitle

\begin{abstract}               
Motivated by neuroscience applications, we introduce the concept of {\it qualitative detection}, that is, the problem of determining on-line the current qualitative dynamical behavior (e.g., resting, oscillating, bursting, spiking etc.) of a nonlinear system. The approach is thought for systems characterized by i) large parameter variability and redundancy, ii) a small number of possible robust, qualitatively different dynamical behaviors and, iii) the presence of sharply different characteristic timescales. These properties are omnipresent in neurosciences and hamper quantitative modeling and fitting of experimental data. As a result, novel control theoretical strategies are needed to face neuroscience challenges like on-line epileptic seizure detection. The proposed approach aims at detecting the current dynamical behavior of the system and whether a qualitative change is likely to occur without quantitatively fitting any model nor asymptotically estimating any parameter. We talk of {\it qualitative detection}. We rely on the qualitative properties of the system dynamics, extracted via singularity and singular perturbation theories, to design low dimensional qualitative detectors.  We introduce this concept on a general class of singularly perturbed systems and then solve the problem for an analytically tractable class of two-dimensional systems with a single unknown sigmoidal nonlinearity and two sharply separated timescales. Numerical results are provided to show the performance of the designed qualitative detector.
\end{abstract}


\end{frontmatter}

\section{Introduction}

In neurosciences, to detect on-line the current activity type of a neuron (spiking, bursting etc.) or of a population of neurons (resting, ictal, interictal, slow/fast oscillations etc.) is of fundamental importance. For epilepsy for instance, such algorithms could be used to detect or even predict seizures. It seems that this problem is very challenging from a control-theoretic viewpoint for several reasons. First, disparate combinations of biophysical parameters are known to lead to the same activity pattern at the cellular level \cite{goldmanetal2001}, and the same degenerated parametrization property propagates at the neuronal circuit level \cite{marder2011}. Second, biophysical parameters change smoothly over time but by doing so they induce sharp, almost discontinuous, transitions between qualitatively different activity types (spiking or bursting, healthy or epileptic, etc.) at the crossing of critical parameter sets. Third, the few available models are often subject to large uncertainties. As a consequence, quantitative modeling and fitting of experimental neuronal data generally constitute an ill-posed problem: a new viewpoint on the problem is needed, which is the purpose of this work.

A first major point is to extract the ruling parameters governing neuronal dynamics and their critical values. In spite of the nonlinear, degenerate structure of the biophysical parameter space, neuronal systems typically exhibit only a few distinct qualitative behaviors. For instance, the activity can be classified as resting, spiking, bursting at the single cell level, and as resting, low amplitude/high frequency, high amplitude/low frequency oscillations at the network level. Predicting the structure of the critical parameter set, where switches in behavior occur, remains an open problem in general. Recently, based on singularity and singular perturbation theories, it was shown that the quantitative parameter space of biophysical neuronal models can be mapped to a small number of lumped parameters. The lumped parameters define low-order polynomial models, capturing the qualitative neuronal behavior at the single cell level \cite{drionetal2015,alessioetalsiam2014}. The resulting geometric framework allows a local, qualitative sensitivity analysis of neuronal dynamics \cite{alessioetalsiam2014,drionetalcdc2015}, which predicts the effect of biophysical parameter variations on the qualitative behavior and the qualitative shape of the critical parameter set. Motivated by the above considerations and inspired by the qualitative sensitivity analysis in \cite{alessioetalsiam2014,drionetalcdc2015}, our first objective is to extract analytically the parameters that rule the activity type of a highly uncertain, redundantly parametrized system and then to detect on-line the activity type and whether a dynamical change is likely to occur.

We formulate these ideas on a general class of nonlinear singularly perturbed control systems that embrace, for instance, the Hodgking-Huxley model for neuronal spiking \cite{HH1952} and the Wilson-Cowan model for neuronal population activity \cite{wilsoncowan1972}. The singularly perturbed nature of the model is justified by neuronal biology, where sharply separated timescales govern the neuronal dynamics~\cite{alessioetalsiam2014,drionetal2015,drionetalcdc2015}. Because the proposed problem is very general and challenging, we then focus on solving it for a two dimensional class of nonlinear systems with a single unknown sigmoidal nonlinearity and two sharply separated timescales. This class of systems captures the qualitative features of the center manifold reduction of higher-dimensional biophysical models, as we explain. In order to cope with the peculiarities of neuronal dynamics, we assume that the exact functional form of the sigmoid nonlinearity is unknown as well as the exact timescale separation. We firstly show that, independently of the particular expression of the sigmoid nonlinearity, the system exhibits either global exponential stability (resting) or relaxation oscillations depending on a single ruling parameter. Guided by singularity theory \cite{golu,alesreal}, we subsequently design an on-line {\it qualitative detector} that  infers the activity type of the system and determines whether a qualitative dynamical transition is close to occur. By qualitative, we mean that the detector neither needs to quantitatively fit any model nor to asymptotically estimate any parameter. Because activity transitions are governed by bifurcations of the underlying vector field, the detection problem put forward here is closely related to the problem of steering a system toward an {\it a priori} unknown bifurcation point \cite{moreauandsontag2003,moreau2003feedback}.

The main contributions of the present work are summarized as follows.\\
$\bullet$Inspired by neurosciences, we formulate the problem of detecting on-line qualitative changes in nonlinear systems for which an exact model is not available.\\
$\bullet$We provide a solution for a class of two-dimensional nonlinear systems, to detect whether the system is close to a transition between resting and oscillatory activity. The stability and the robustness of the detector are ensured via Lyapunov techniques and geometric singular perturbation arguments.\\
$\bullet$We provide numerical evidences that the same detector performs well when tested on a classical, high dimensional, biophysical neuronal model of resting and spiking oscillations due to Hodgkin and Huxley \cite{HH1952}.

The paper is organized as follows. In Section~\ref{sec:2}, we first formulate the idea on a general nonlinear systems. We then introduce and analyze  a class of singularly perturbed nonlinear systems with a single unknown sigmoidal nonlinearity. In Section~\ref{sec:3}, we design and analyze the qualitative detector for the latter. Numerical simulations are presented in Section \ref{sec:4}. Section \ref{sec:5} concludes the paper. The proofs are given in the appendix.

\textbf{Notation.} Let $\RR_{\geqslant 0}:=[0,\infty)$, $\RR_{>0}:=(0,\infty)$, and $\mathbb N:=\{1,2,\ldots,\}$. The usual Euclidean norm is denoted by $\vert \cdot\vert$. For $x,y,z\in\RR$, the vector $[x,y,z]^\top$ is denoted by $(x,y,z)$. For a function $h:\RR_{\geqslant0}\to\RR^n$, the associated infinity norm is denoted by $\Vert h\Vert_{[0,\infty)}=\sup_{s\in[0,\infty)}\vert h(s)\vert$, when it is well-defined. We use $\text{sgn}(x)$ to denote the sign function from $\RR$ to $\{-1,0,1\}$ with $\text{sgn}(0)=0$. For any function $f:\RR \to \RR$, we denote the range of $f$ as $f(\RR)=\{z: z=f(x) \;\text{for some}\; x\in \RR\}$. Let $A, B$ be two non-empty subsets in $\RR^n$, their Hausdorff distance is noted by $d_H(A,B)=\max\{\sup_{a\in A}\inf_{b\in B} \vert a-b\vert,\;\sup_{b\in B}\inf_{a\in A} \vert a-b\vert\}$. The point to set distance from $p\in\RR^n$ to $A\subseteq \RR^{n}$ is denoted by $\vert p\vert_A=\inf\limits_{a\in A}\vert p-a\vert$. 
\section{Problem statement}\label{sec:2}
\subsection{General formulation}\label{sec21}
We first introduce the main ideas on the following general class of nonlinear singularly perturbed systems
\begin{subequations}\label{EQ: generic nonlinear}
\begin{align}[left = \!\!\!\!\empheqlbrace\,]
\label{gn1}
\dot x_f&=f(x_f,h_s(x_s),\alpha,u),\\\label{gn2}
\dot x_s&=\varepsilon g(h_f(x_f),x_s),\\\label{gn3}
y&=(h_f(x_f),h_s(x_s)),
\end{align}
\end{subequations}
where $x_f\in\mathbb{R}^n$, $x_s\in\mathbb{R}^m$ are the state variables, $\alpha\in\mathbb{R}^p$ is an {\it unknown} parameter vector, $u\in\mathbb{R}$ is a constant input, $0<\varepsilon\ll 1$ is an {\it unknown} parameter, $h_s:\mathbb{R}^m\to\mathbb{R}$, $h_f:\mathbb{R}^n\to\mathbb{R}$ and all functions are smooth.

System~\eqref{EQ: generic nonlinear} is said to be {\it redundantly parametrized} if it depends on numerous parameters, i.e. $p\gg1$, but the possible qualitatively distinct dynamical behaviors are only a few. Parametric redundancy naturally arises under two conditions: timescale separation and the presence of an organizing center, as explained below.

System (\ref{EQ: generic nonlinear}) evolves  according to two time scales, since $\varepsilon$ is small. Under suitable stability and monotonicity conditions of the fast dynamics (see e.g. \cite{gedeonsontag2007}), the trajectories of the fast variable $x_f$ are forced toward its instantaneous quasi-steady state $x_f^*(t)$, which is defined as follows. For all $t\in\mathbb R_{\geq0}$, $\alpha\in\mathbb R^p$, $u\in\mathbb R$, let $x_s(t)$ be the slow state at time $t$. Then $x_f^*(t)$ is a solution of the quasi-steady state equation
$f(x_f^*(t),h_s(x_s(t)),\alpha,u)=0$,
and the fast output $y_f=h_f(x_f)$ will change accordingly.
Now, for $u$ fixed, define the map $G:\mathbb{R}^{n+1}\times\mathbb{R}\times\mathbb{R}^{p}\to\mathbb{R}^{n+1}$ as
\begin{equation*}
\label{unfolding}
G(x,\lambda,\alpha):=\begin{bmatrix}
f(x_f,\lambda,\alpha,u)\\
y_f-h_f(x_f)
\end{bmatrix},
\end{equation*}
where $x=(x_f,y_f)$ and $\lambda$ is the {\it bifurcation parameter}. For system \eqref{EQ: generic nonlinear} to be redundantly parameterized, we assume that the map $G$ is {\it organized by a $k$-codimensional singularity}, with $k\ll p$, in the sense of \cite[Section III.1]{golu}. This means that the qualitative shape of the graph of the, possibly multivalued, mapping $\lambda\mapsto x_f$ defined by solving the equation $G(x,\lambda,\alpha)=0$ is given by one of the {\it persistent bifurcation diagrams} \cite[Sections III.5,6]{golu} in the universal unfolding of the organizing singularity. This shape is exactly what determines the quasi-steady state evolution of the fast variable. Under suitable stability and monotonicity properties of the slow dynamics, (for instance that they provide slow adaptation - negative feedback - on the fast dynamics), it follows that the possible qualitative dynamic output behaviors of system~\eqref{EQ: generic nonlinear} are fully determined by the persistent bifurcation diagrams of the organizing singularity. Hence, since $k\ll p$, the possible different persistent bifurcation diagrams, and therefore the possible output behaviors are much less than the number of model parameters: parametric redundancy occurs.

When system (\ref{EQ: generic nonlinear}) exhibits parametric redundancy, we do not need to know the full vector of parameters $\alpha\in\mathbb R^p$ to detect what is the current dynamic behaviour of the system and whether a qualitative change is prone to appear. The projection of  $\alpha$ onto the $k$-dimensional unfolding space of the organizing suffices to determine the model behavior. The sensitivity analysis proposed in \cite{drionetal2015} is fully based on this idea. In this paper, we explore this viewpoint to formulate the following {\it qualitative detection} problem:
Let $f$, $g$, $h_f$, $h_s$ be unknown. Let the organizing singularity be known. Let the slow and fast output be measurable. Let the dimension $p$ of the parameter space be {\it unknown}. Can we detect on-line in which region of the organizing singularity unfolding space the parameter vector lies and whether this vector is approaching a transition variety, where a qualitative dynamical transition occurs?

This problem is new, general and very challenging. That is the reason why we solve it for a particular yet relevant class of systems of the form of (\ref{EQ: generic nonlinear}) in the following.

\subsection{A tractable class of systems}
In the following, we concentrate on systems of the form 
\begin{subequations}\label{EQ: sigmoidal hysteresis class}
\begin{align}[left = \!\!\!\!\empheqlbrace\,]
\label{1_1}
\dot x_f&=-x_f+S(\beta x_f+u-x_s),\\\label{1_2}
\dot x_s&=\varepsilon(x_f-x_s),
\end{align}
\end{subequations}
where $x_f, x_s\in\mathbb R$ are the state variables, $\beta\in\RR$ is an {\it unknown} parameter, $u\in \RR$ is the input, which is assumed to be known and constant, and $0<\varepsilon\ll1$ is small {\it unknown} parameter. The mapping $S:\RR \to \RR$ is an {\it unknown} sigmoid function, which is assumed to satisfy the following properties. 
\begin{assumption}
\label{sigmo_assum}
The following properties holds.\\
a$)$ $S$ is smooth on $\RR$.\\
b$)$ $S(0)=0$.\\
c$)$ $S^\prime (a)>0$ for all $a\in\mathbb R$ (monotonicity), and ${\rm argmax}_{a\in\mathbb R} S^\prime(a)=0$ (sector-valued).\\
d$)$ ${\rm sgn}\left(S^{\prime\prime}(a)\right)=-{\rm sgn}(a)$ for all $a\in\mathbb R$.$\hfill \Box$
\end{assumption}
Standard sigmoid functions such as $a\mapsto c_1 \tanh(c_2 a)$, $a\mapsto \frac{c_1}{1+e^{-c_2 a}}-\frac{c_1}{2}$, $a\mapsto e^{-e^{-c_1a}}-c_2$ with $c_1,c_2\in\RR_{>0}$, verify Assumption \ref{sigmo_assum}. Moreover, when $S=\sum\limits^{n}_{i=1}S_i$ where $S_1,\cdots, S_i$ with ($i\in  \mathbb N$) are sigmoid functions satisfying Assumption \ref{sigmo_assum}, then so does $S$.

Due to the small parameter $\varepsilon$, system (\ref{EQ: sigmoidal hysteresis class}) evolves according to two time scales. We follow the standard approach of singular perturbation theory \cite{SPcontrol, Fenichel1979a} to decompose system (\ref{EQ: sigmoidal hysteresis class}) into two subsystems, which represent the fast dynamics and the slow dynamics, respectively, called the {\it layer} and the {\it reduced} subsystems. By setting $\varepsilon=0$ in (\ref{1_2}), we obtain the layer dynamics
\begin{subequations}\label{EQ: sigmoidal hysteresis class layer}
\begin{align}[left = \empheqlbrace\,]\label{lay_1}
\dot x_f&=-x_f+S(\beta x_f+u-x_s),\\\label{lay_2}
\dot x_s&=0.
\end{align}
\end{subequations}
Because of (\ref{lay_2}), the variable $x_s$ can be treated as a parameter in (\ref{lay_1}). 

To account for the slow variations of $x_s$, we rescale the time as $\tau=\varepsilon t$, hence $\frac{d}{dt}=\varepsilon\frac{d}{d\tau}$. Then, system (\ref{EQ: sigmoidal hysteresis class}) becomes
\begin{subequations}\label{trans_red}
\begin{align}[left = \empheqlbrace\,]
\label{trans_1}
\varepsilon x^\prime_f&=-x_f+S(\beta x_f+u-x_s),\\\label{trans_2}
x^\prime_s&=x_f-x_s,
\end{align}
\end{subequations}
where $\prime$ stands for $\frac{d}{d\tau}$. Setting $\varepsilon=0$ in (\ref{trans_1}), we obtain the reduced dynamics
\begin{subequations}\label{EQ: sigmoidal hysteresis class reduced}
\begin{align}[left = \empheqlbrace\,]
\label{rec_1}
0&=-x_f+S(\beta x_f+u-x_s)\\\nonumber
&=:F_f(x_f,u-x_s,\beta),\\\label{rec_2}
x^\prime_s&=x_f-x_s=:F_s(x_f,x_s).
\end{align}
\end{subequations}
The reduced dynamics evolves in the {\it slow time} $\tau$ and is an algebro-differential equation. In particular it defines a one-dimensional vector field over the {\it critical manifold}
\begin{IEEEeqnarray}{rCl}\label{EQ: sigmoidal hysteresis crit man}
\IEEEyesnumber
S^0_{u,\beta}&:=&\{(x_s,x_f)\in\mathbb{R}^2:\ F_f(x_f,u-x_s,\beta)=0\}.
\end{IEEEeqnarray}
The critical manifold $S_{u,\beta}^0$ depends on $u$ and $\beta$. However, since we assume $u$ and $\beta$ constant, for simplicity we omit the index $u$ and $\beta$ in the rest of the paper. The critical manifold will therefore simply be denoted by $S^0$.
The idea behind singular perturbations is that dynamics (\ref{EQ: sigmoidal hysteresis class}) are well approximated by layer dynamics far from $S^0$ and by trajectories of the reduced dynamics close to $S^0$. The shape of $S^0$ therefore plays a key role in determining the behavior of system (\ref{EQ: sigmoidal hysteresis class}). 

\begin{figure}[H]
\begin{subfigure}[t]{0.47\columnwidth}
\centering
\includegraphics[width=\columnwidth,height=3cm]{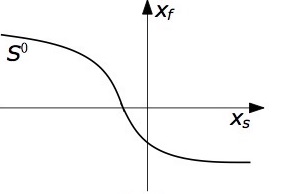}
\caption{$\beta<\beta_c$}\label{simuxffix}
\end{subfigure}
\quad
\begin{subfigure}[t]{0.47\columnwidth}
\centering
\includegraphics[width=\columnwidth,height=3cm]{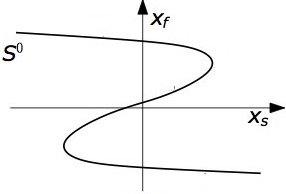}
\caption{$\beta>\beta_c$}\label{simuxflc}
\end{subfigure}
\caption{Qualitative shape of the critical manifold $S^0$.}\label{FIG: hy crit manifold}
\end{figure}

Let 
\begin{equation}
\label{bc}
\beta_c:=\frac{1}{S^\prime(0)},
\end{equation}
which is well-defined and strictly positive according to item c$)$ of Assumption \ref{sigmo_assum}. The parameter $\beta_c$ is crucial in shaping the critical manifold $S^0$, as illustrated in Figure \ref{FIG: hy crit manifold}. For $\beta<\beta_c$ the shape of $S^0$ is given by graph of a monotone decreasing function. Whereas for $\beta>\beta_c$ the shape of $S^0$ is $S$-shaped. We emphasize that $\beta_c$ is unknown since so is the sigmoid function $S$. 

The qualitative properties of the critical manifold fully determine the model dynamical behaviors, as formalized in the next proposition, whose proof can be found in \cite{qualitativefull2}.

\begin{proposition}\label{LEM: stable osci dichotomy}\
Consider system (\ref{EQ: sigmoidal hysteresis class}), the following holds.\\
{\rm 1)} For any $\beta-\beta_c<0$ and any constant input $u$, there exists $\bar\varepsilon>0$ such that, for all $\varepsilon\in(0,\bar\varepsilon]$, system (\ref{EQ: sigmoidal hysteresis class}) has a globally exponentially stable fixed point. Moreover, all trajectories converge in an $\mathcal O(\varepsilon)$-time to an $\mathcal O(\varepsilon)$-neighborhood of the critical manifold $S^0$ defined in (\ref{EQ: sigmoidal hysteresis crit man}).\\
{\rm 2.a)} For all $0<\beta-\beta_c<1$, there exists $\bar u>0$, such that, for any constant input $u$ in $(-\bar u,\bar u)$, there exists $\bar\varepsilon>0$ such that, for all $\varepsilon\in(0,\bar\varepsilon]$, system (\ref{EQ: sigmoidal hysteresis class}) possesses an almost globally asymptotically stable{\footnote{We adopt the notion of almost global asymptotical stability defined in \cite{angeli2004}}} and locally exponentially stable periodic orbit $P^\varepsilon$.\\
{\rm 2.b)} Let $\beta, \bar u, u, \bar \varepsilon, \varepsilon$ be as in item {\rm2.a)}, and let $T^\varepsilon$ be the period of $P^\varepsilon$ and $p:\mathbb R\to P^\varepsilon$, with $p(t)=p(t+T^\varepsilon)$ for all $t\in\mathbb R$, be a periodic solution of (\ref{EQ: sigmoidal hysteresis class}). Let $P^0$ be the singular limit of $P^\varepsilon$. The Hausdorff distance between $P^\varepsilon$ and $P^0$ verifies $d_H(P^\varepsilon,P^0)=\mathcal O(\varepsilon^{2/3})$. Moreover, for all $t_0\in\mathbb R$, there exists $\delta T^\varepsilon\subset [t_0,t_0+T^\varepsilon)$ and $|\delta T^\varepsilon|= \mathcal O(\varepsilon) T^\varepsilon$, such that $|p(t)-S^0|<\mathcal O(\varepsilon)$ for all $t\in[t_0,t_0+T^\varepsilon)\setminus \delta T^\varepsilon $. $\hfill \Box$
\end{proposition}
Proposition \ref{LEM: stable osci dichotomy} shows that system (\ref{EQ: sigmoidal hysteresis class}) either has a globally exponentially stable fixed point (item 1)), which we call the {\it resting activity}, or almost all its solutions converge to a limit cycle (item 2)), in which case we talk of {\it oscillating activity}. The nature of the activity is solely determined by the sign of $\beta-\beta_c$, independently of the particular expression of $S$.

\begin{remark}
Proposition \ref{LEM: stable osci dichotomy} states stability results for $\beta<\beta_c+1$. When $\beta>\beta_c+1$, system (\ref{LEM: stable osci dichotomy}) exhibits three equilibria, two of which are asymptotically stable and one is unstable. Then trajectories converge to the two stable equilibria. \end{remark}

\begin{remark}
\label{remHH}
Higher-dimensional control systems of the form (\ref{EQ: generic nonlinear}) may be reduced to the simpler lower-dimensional system (\ref{EQ: sigmoidal hysteresis class}) via center manifold reduction \cite{Guckenheimer1983}. We sketch the ideas as follows. The transition from globally asymptotically stable fixed point to almost globally asymptotically, locally exponentially stable oscillations in model (\ref{EQ: sigmoidal hysteresis class}) is via a Hopf bifurcation, as it can easily be verified by linearization at the origin (for $u=0$). Assume that the same bifurcation exists in the full model (\ref{EQ: generic nonlinear}). Then, close to this bifurcation, by the center manifold theorem, the dynamics of model (\ref{EQ: generic nonlinear}) can be reduced to a two-dimensional invariant manifold: the center manifold of the Hopf bifurcation. The same assumptions as for (\ref{EQ: sigmoidal hysteresis class}) can then be formulated for the full model (\ref{EQ: generic nonlinear}) by restricting them to the Hopf center manifold. Classical two-dimensional reduction of the Hodgkin-Huxley model \cite{fitzhugh1961,rinzel1985} satisfy these conditions with the sodium maximal conductance $\bar g_{Na}$ as unfolding parameter and the applied current $I_{app}$ as constant input. To illustrate this, numerical results for Hodgkin-Huxley model are provided in Section \ref{sec4.2}\end{remark}

In the next section we reformulate the qualitative detection problem in the specific context of model~\eqref{EQ: sigmoidal hysteresis class} and we provide a solution.

\section{Qualitative detector}\label{sec:3}
\subsection{Objective and design}
Our objective is to detect on-line the current type of activity of system (\ref{EQ: sigmoidal hysteresis class}), i.e. whether the solution converges to a fixed point or to a limit cycle, and how close the system is to a change of activity. In view of Proposition~\ref{LEM: stable osci dichotomy}, activity switch in system (\ref{EQ: sigmoidal hysteresis class}) is ruled by $\beta$ only. Hence, the most natural idea is to estimate $\beta$. For this purpose, we might approximate the sigmoid function $S$ by a Taylor expansion up to some order and quantitatively estimate the expansion coefficients together with $\varepsilon$ and $\beta$. However, due to the generic impossibility of quantitatively modeling and fitting experimental neuronal data, this approach would not overcome, in a real experimental setting, the limitation of quantitative/asymptotic estimation. Instead, we aim at {\it qualitatively} detection the different dynamical activities of the system: whether the system exhibits oscillation or has a globally exponentially stable fixed point and whether it is near to a change of activity.

Proposition \ref{LEM: stable osci dichotomy} proves that (almost) all trajectories of system (\ref{EQ: sigmoidal hysteresis class}) converge to an $\mathcal O(\varepsilon)$-neighborhood of the critical manifold  $S^0$ in both the fixed point and the limit cycle cases, at least for most of the time in the latter case. Based on singularity theory \cite{golu,alesreal} and in view of Assumption~\ref{sigmo_assum}, the shape of the critical manifold $S^0$ is, modulo an equivalence transformation, the same as that of the set
\begin{equation}
\label{eq_qualitative}
\{(x_f,x_s):-x_f^3+(\beta-\beta_c)x_f+u-x_s=0\},
\end{equation}
whose definition is independent of $S$. We exploit this information to design and analyze our detector in the following.

We first make the following assumption.
\begin{assumption}
\label{est_assm1}
Both $x_f$ and $x_s$ variables are measured in system (\ref{EQ: sigmoidal hysteresis class}).
\end{assumption}

\begin{remark}
When only the fast dynamics is available by measurement, we may reconstruct the slow variable using a state estimator of the form $\dot {\hat{x}}_s=\hat \varepsilon (x_f-\hat{x}_s)$, for some small $\hat{\varepsilon}>0$. The analysis of this scenario is not addressed in this paper and is left for future work. 
\end{remark}

We propose the nonlinear qualitative detector
\begin{equation}
\label{est_xi}
\!\!\!\!\!\dot{\hat\beta}=\!-kx_f(-x_f^3+\hat\beta x_f+u-x_s)\!=:\hat f(\hat\beta,x_f,u-x_s),
\end{equation}
where $k>0$ is a design parameter. The definition of $\hat f$ does not require the knowledge of $S$ or $\varepsilon$, it only relies on the information  provided by (\ref{eq_qualitative}).

In the following, we define $\hat\beta^*$ the steady state of $\hat\beta$, and we analyze its relationship with respect to $\beta-\beta_c$. We use the terminology steady-state with some slight abuse as we will see that $\hat\beta^*$ is not always a constant. We then investigate the incremental stability properties of (\ref{est_xi}), which we finally exploit to prove the stability of (\ref{est_xi}) and the convergence of $\hat\beta$ to $\hat\beta^*$.

\subsection{Steady state properties of the qualitative detector}\label{sec:32}
Steady state of $\hat\beta$ with $x_f\neq 0$ satisfies $-x_f^3+\hat\beta x_f+u-x_s=0$ in view of (\ref{est_xi}). We thus implicitly define the steady state function $\hat\beta^*(x_f,x_s,u)$ by
\begin{equation}
\label{betahatstar1}
-x^3 _f+\hat\beta^*(x_f,x_s,u)x_f+u-x_s=0,
\end{equation}
for $u,\beta \in\RR$ and $ x_f, x_s\in S^0$, i.e. when $x_f, x_s$ are related via the critical manifold equation
\begin{equation}
\label{betahatstarcrit}
-x_f+S(\beta x_f+u-x_s)=0.
\end{equation}
Recalling that the sigmoid function $S$ is invertible on $S(\RR)$ according to item c$)$ of Assumption \ref{sigmo_assum}, we can explicitly invert (\ref{betahatstarcrit}) as follows
\begin{equation}
\label{sol_xs}
x_s=-S^{-1}(x_f)+\beta x_f+u.
\end{equation}
Replacing this expression for $x_s$ in (\ref{betahatstar1}), we obtain the explicit expression of $\hat\beta^*(x_f,x_s,u)$
\begin{equation}
\label{ift_sin2}
\hat\beta^*(x_f,x_s,u) =\frac{x^{3}_f-S^{-1}(x_f)+\beta x_f}{x_f}.
\end{equation}
Because $\hat\beta^*$ depends only on $x_f$ and $\beta$ in (\ref{ift_sin2}), we write $\hat\beta^*(x_f,\beta)$ in the sequel.

The next theorem provides important properties of $\hat\beta^*(x_f,\beta)$. In particular, it ensures that $\hat \beta^*(x_f,\beta)$ is well-defined when $x_f=0$. Furthermore, it formally characterises in which sense $\hat\beta^*$ provides qualitative information about $\beta-\beta_c$.
\begin{theorem}
\label{property_quali_betahs}
The following holds.\\
{\rm1)} $\hat\beta^*$ is smooth on $S(\RR)\times \RR$.\\
{\rm2)} For all $\beta<\beta_c$ and $x_f\in S(\RR)$,
$\frac{\partial\hat\beta^*(x_f,\beta)}{\partial\beta}>0$.\\
{\rm3)} For all $\beta\in(0,\beta_c+1)$ and $x_f$ in a neighborhood of the origin, $\frac{\partial\hat\beta^*(x_f,\beta)}{\partial\beta}>0$. Moreover, let $x^+_{fold}>0$ and $x^-_{fold}<0$ be the two distinct solutions to $-1+\beta S^\prime(S^{-1}(x_{f}))=0$. If, for $x_{fold}\in\{x^+_{fold}, x^-_{fold}\}$,
\begin{equation}
\label{cond_suf}
2x^2_{fold}+\frac{S^{-1}(x_{fold})}{x_{fold}}-(S^{-1})^\prime (x_{fold})>0, 
\end{equation}
then $\frac{\partial \hat\beta^*(x_f,\beta)}{\partial \beta}>0$ for $x_f\in(-\infty,x^-_{fold})\cup (x^+_{fold},+\infty)$.\\
{\rm4)} $\hat\beta^*(x_f,\beta)=\beta-\beta_c+\mathcal O(x^2_f)$ for any $x_f\in S(\RR)$ and $\beta\in \RR$.$\hfill \Box$
\end{theorem}

Theorem~\ref{property_quali_betahs} reveals that $\hat\beta^*$ provides information about the qualitative trend of $\beta$. According to item 2$)$ of Theorem~\ref{property_quali_betahs}, an increase in $\hat\beta^*$ indicates that $\beta$ is increasing in the resting activity (i.e. $\beta<\beta_c$). If a change of activity is prone to appear (i.e. $\beta$ is close to $\beta_c$), then necessarily $x_f$ will lie in a neighborhood of the origin (for $u$ small), which relates to the point of infinite slope, corresponding to the hysteresis singularity of the critical manifold. In this case, $\hat\beta^*$ will change its sign from negative to positive by virtue of item 4$)$ of Theorem~\ref{property_quali_betahs}, and then it will keep increasing in $\beta$ in the oscillation activity (i.e. $\beta_c<\beta<\beta_c+1$) in view of item 3$)$ of Theorem~\ref{property_quali_betahs}. The same holds when $\hat\beta^*$ goes from positive to negative values. Hence, properties of $\hat\beta^*$ ensure that we can detect online when $\beta$ is about to cross, either from above or from below, the critical value $\beta_c$, at least for small $u$. Condition (\ref{cond_suf}) can be verified off line for different sigmoid functions and might be instrumental in tracking variations of $\beta$ for $\beta-\beta_c$ strictly positive and large.

\subsection{Incremental stability properties of the detector}\label{sec:33}
If the value of $\hat\beta$ converges to $\hat\beta^*$ as time tends to infinity, then (\ref{est_xi}) can inform us about the activity of system (\ref{EQ: sigmoidal hysteresis class}) as explained above. To prove this, we first show that system (\ref{est_xi}) exhibits incremental stability properties \cite{angeli2002} if the input signal $x_f$ to the detector is bounded and a persistency of excitation (PE) condition holds. We then prove that under mild conditions, the solutions to system (\ref{EQ: sigmoidal hysteresis class}) do satisfy the aforementioned PE requirement.

We use the following notion of PE, which is similar to the one in \cite{anderson1977}. 

\begin{definition}
The piecewise continuous signal $y:\RR_{\geqslant 0}\rightarrow \RR$ is $(T,\mu)$-PE, if there exist $T, \mu>0$ such that $ \int_t^{t+T} y^2(\tau) \;d\tau\geqslant \mu$ for all $t\geqslant 0$.$\hfill \Box$
\end{definition}

\begin{proposition}
\label{proposition2}
For any $\Delta,M, T,\mu>0$, there exist $\ell_1(\Delta,M,T,\mu),\ell_2(\Delta,M,T,\mu), \ell_3(\Delta,M,T,\mu)>0$ such that, if (i) $x_{f,i}$ is $(T,\mu)$-PE, $i\in \{1,2\}$, and $(ii)$ $x_{f,i},x_{s,i},u_i$ are piecewise continuous with \\$\max(\Vert x_{f,i}\Vert_\infty,\Vert x_{s,i}\Vert_\infty,\Vert u_i\Vert_\infty)\leqslant M$, $i\in \{1,2\}$, then for any $\vert\hat\beta_i(0)\vert\leqslant \Delta$, $i\in \{1,2\}$, the corresponding solution $\hat\beta_i$ to system (\ref{est_xi}) satisfies for all $t\geqslant 0$
\begin{eqnarray}
\nonumber
\vert \hat\beta_1(t)-\hat\beta_2(t)\vert &\leqslant& \ell_1(\Delta,M,T,\mu)e^{-\ell_2(\Delta,M,T,\mu)t} \vert \hat\beta_1(0)-\hat\beta_2(0)\vert\\\nonumber
&+&\ell_3(\Delta,M,T,\mu) \biggl(\Vert x_{f1}-x_{f2}\Vert_{[0,t)}\\\label{diss_eq}
&+&\Vert x_{s1}-x_{s2}\Vert_{[0,t)}+\Vert u_{1}-u_{2}\Vert_{[0,t)}\biggr).
\end{eqnarray}$\hfill\Box$
\end{proposition}
Proposition \ref{proposition2} states that system (\ref{est_xi}) satisfies a semiglobal input-to-state incremental stability property, when its input signals verify conditions $(i)$ and $(ii)$. In the next lemma, we analyze under what conditions the solution to system (\ref{EQ: sigmoidal hysteresis class}) verifies these two conditions.
\begin{lemma}
\label{lem_PE}
Consider system (\ref{EQ: sigmoidal hysteresis class}).\\
{\rm1)} For any $\beta-\beta_c<0$, any $\Delta>\delta>0$ and any constant input u satisfying $\delta<\vert u\vert<\Delta$, there exist $\bar\varepsilon, T(\Delta,\delta), \mu(\Delta,\delta), C(\Delta,\delta)>0$, such that for all $\varepsilon\in(0,\bar\varepsilon]$ and initial conditions $x_f(0)$ and $x_s(0)$ with $\vert x_f(0)\vert, \vert x_s(0)\vert\leqslant \Delta$, the corresponding solution to system (\ref{EQ: sigmoidal hysteresis class}) is such that $(i)$ its $x_f$-component is $(T(\Delta,\delta),\mu(\Delta,\delta))$-PE and\\ $(ii)$ $\max(\Vert x_f\Vert_{[0,\infty)},\Vert x_s\Vert_{[0,\infty)},\Vert u\Vert_{[0,\infty)})\leqslant C(\Delta,\delta)$.\\
{\rm2)} For any $0<\beta-\beta_c<1$, any $\Delta>\delta>0$ there exists $\bar u\in(0,\Delta)$ such that for any $u\in(-\bar u,\bar u)$, there exist $\bar\varepsilon,T(\Delta,\delta), \mu(\Delta,\delta), C(\Delta,\delta)>0$, such that for all $\varepsilon\in(0,\bar\varepsilon]$ and initial conditions $x_f(0)$ and $x_s(0)$ with $\vert x_f(0)\vert, \vert x_s(0)\vert\leqslant \Delta$ and $\vert x_f(0)-p^*_f\vert\geqslant \delta, \vert x_s(0)-p^*_s\vert\geqslant \delta$, where $(p^*_f,p^*_s)$ is the unique unstable fixed point of (\ref{EQ: sigmoidal hysteresis class}), the corresponding solution to system (\ref{EQ: sigmoidal hysteresis class}) satisfies conditions $(i)$ and $(ii)$ as in the above item 1). $\hfill \Box$
\end{lemma}

Lemma \ref{lem_PE} means that the conditions of Proposition \ref{proposition2} are verified semiglobally by solutions to system (\ref{EQ: sigmoidal hysteresis class}) whenever $\beta<\beta_c+1$, under mild conditions on the input $u$.
\subsection{Main result}
To finalize the analysis, we consider the overall system, which is a cascade of system (\ref{EQ: sigmoidal hysteresis class}) with (\ref{est_xi})
\begin{subequations}\label{3D}
\begin{align}[left = \!\!\!\!\empheqlbrace\,]
\label{3D1_1}
\dot x_f&=-x_f+S(\beta x_f+u-x_s),\\\label{3D1_2}
\dot x_s&=\varepsilon(x_f-x_s),\\\label{3D1_3}
\dot {\hat\beta}&=-kx_f(-x_f^3+\hat\beta x_f+u-x_s).
\end{align}
\end{subequations}
We state the stability property of system (\ref{3D}) and the relationship between $\hat\beta$ and $\hat\beta^*$ in the following two theorems respectively.

\begin{theorem}
\label{thm_1}
Consider system (\ref{3D}).
\\
{\rm 1)} For any $\beta<\beta_c$, $\Delta>\delta>0$, and any constant input $u$ satisfying $\delta< \vert u\vert<\Delta$, there exists $\bar\varepsilon>0$ such that, for any $\varepsilon\in(0,\bar\varepsilon]$, system (15) has a globally asymptotically stable and locally exponentially stable fixed point.
\\
{\rm 2)} For any $0<\beta-\beta_c<1$, $\Delta>0$, there exist $\bar\varepsilon>0$ and $\bar u\in(0,\Delta)$ such that for any $u\in(-\bar u,\bar u)$, and for any $\varepsilon\in(0,\bar\varepsilon]$, system (15) possesses an almost globally asymptotically stable and locally exponentially stable periodic orbit $Q^\varepsilon$ with period $T^\varepsilon$. $\hfill \Box$
\end{theorem}

Theorem \ref{thm_1} shows that the solution to the overall system (\ref{3D}) converges either to a stable fixed point or, almost all of these, converge to a limit cycle depending on the value of $\beta$. In the next theorem, we specify the relationship between $\hat\beta$ and $\hat\beta^*$.
\begin{theorem}
\label{thm_betah}
Consider system (\ref{3D}). \\
{\rm 1)} For any $\beta<\beta_c$, $\Delta>\delta>0$, and any constant input $u$ satisfying $\delta< \vert u\vert<\Delta$, there exists $\bar\varepsilon>0$ such that, for any $\varepsilon\in(0,\bar\varepsilon]$ and any initial conditions $\vert x_f(0)\vert, \vert x_s(0)\vert,\vert \hat\beta(0)\vert\leqslant \Delta$, $\hat\beta(t)$ converges to $\hat\beta^*(x_f^*,\beta)$ as time goes to infinity, where $x^*_f$ is $x_f$-component of the globally asymptotically stable and locally exponentially stable fixed point of (\ref{3D}), as stated in item 1) of Theorem \ref{thm_1}.\\
{\rm 2)} For any $0<\beta-\beta_c<1$, $\Delta>\delta>0$, there exist $\bar\varepsilon>0$ and $\bar u\in(0,\Delta)$ such that for any $u\in(-\bar u,\bar u)$, any $\varepsilon\in(0,\bar\varepsilon]$ and any initial conditions $\vert x_f(0)\vert, \vert x_s(0)\vert, \vert \hat\beta(0)\vert\leqslant \Delta$ and $\vert x_f(0)-p^*_f\vert\geqslant \delta, \vert x_s(0)-p^*_s\vert\geqslant \delta$, where $(p^*_f,p^*_s)$ is the unique unstable fixed point of system (\ref{3D1_1})-(\ref{3D1_2}), $\hat\beta(t)$ converges in an $\mathcal O(\varepsilon)$-neighborhood of $\hat\beta^*(x_f^{lc}(t+\theta),\beta)$ for all $t\in[t_0,t_0+T^\varepsilon)\setminus\delta T^\varepsilon $ with $|\delta T^\varepsilon|= \mathcal O(\varepsilon) T^\varepsilon$, where $x_f^{lc}(t+\theta)$ is $x_f$-component of the $T^\varepsilon$-periodic almost globally asymptotically stable and locally exponentially stable periodic orbit of (\ref{3D}), as stated in item 2) of Theorem \ref{thm_1} with phase $\theta\in[0,2\pi)$.
$\hfill \Box$
\end{theorem}
When $\beta<\beta_c$, item \rm{1)} of Theorem \ref{thm_betah} shows that the value of $\hat\beta$ is $\hat\beta^*$ when time goes to infinity. For $0<\beta-\beta_c<1$, item \rm{2)} of Theorem \ref{thm_betah} indicates that $\hat\beta$ is in an $\mathcal O(\varepsilon)$-neighborhood of $\hat\beta^*$ for all the time, after a sufficiently long time, except during jumps of length $|\delta T^\varepsilon|= \mathcal O(\varepsilon) T^\varepsilon$. We see that $\hat\beta^*$ is not a constant in this case but a state-dependent signal. This is not an issue as all what matters is the sign of $\hat\beta^*$, which is informative according to Theorem \ref{property_quali_betahs}, and not its actual value. Various filtering techniques, like simple averaging over a sliding window, can be used to extract an approximated constant value of $\hat\beta^*$ if needed in this case, see Section \ref{sec4.1} for instance.

\begin{remark}
\label{rubst}
In view of \cite{sontagwang1996}, the stability property of system (\ref{3D}) in Theorem \ref{thm_1} implies its local input-to-state stability with respect to any small exogenous disturbance.  In other words, the detector is robust to small measurement noises as well as small perturbations in its dynamics.
\end{remark}

\section{Numerical illustrations}\label{sec:4}

\subsection{System (\ref{EQ: sigmoidal hysteresis class})}\label{sec4.1}
\begin{figure}[H]
\centering
\begin{subfigure}[t]{0.47\columnwidth}
\centering
\includegraphics[width=\columnwidth,height=3.5cm]{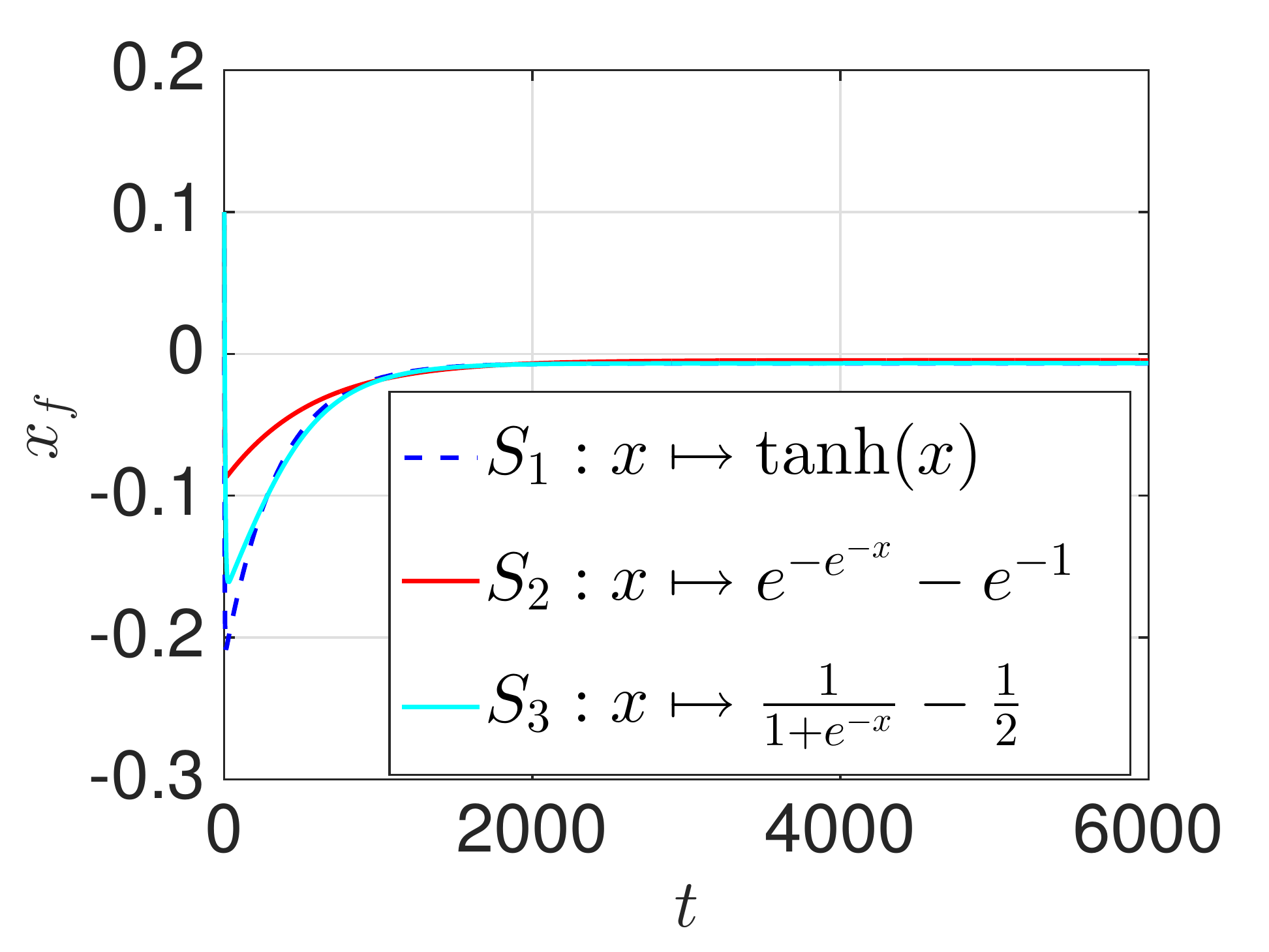}
\end{subfigure}
\quad
\begin{subfigure}[t]{0.47\columnwidth}
\centering
\includegraphics[width=\columnwidth,height=3.5cm]{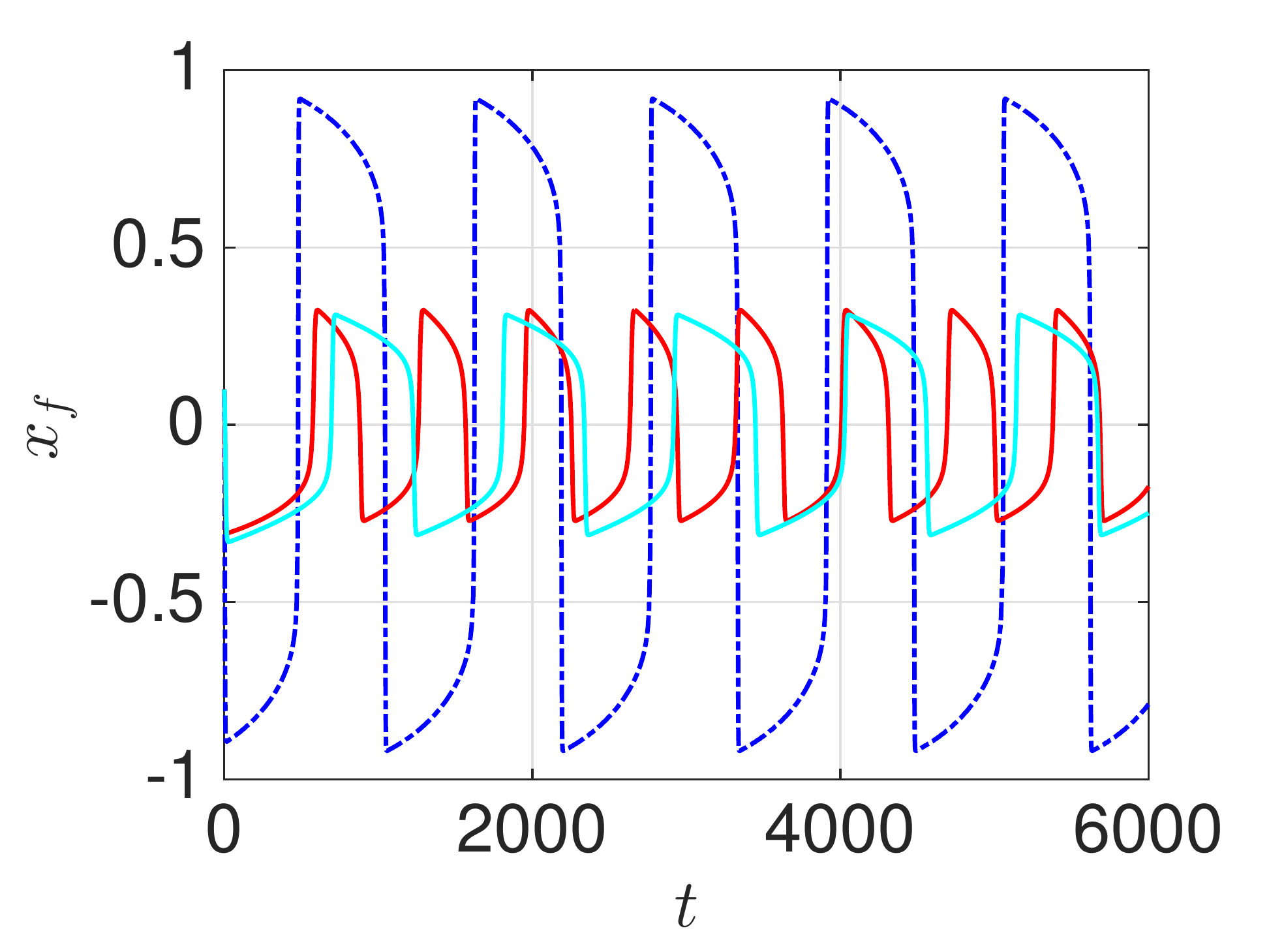}
\end{subfigure}
\begin{subfigure}[t]{0.47\columnwidth}
\centering
\includegraphics[width=\columnwidth,height=3.5cm]{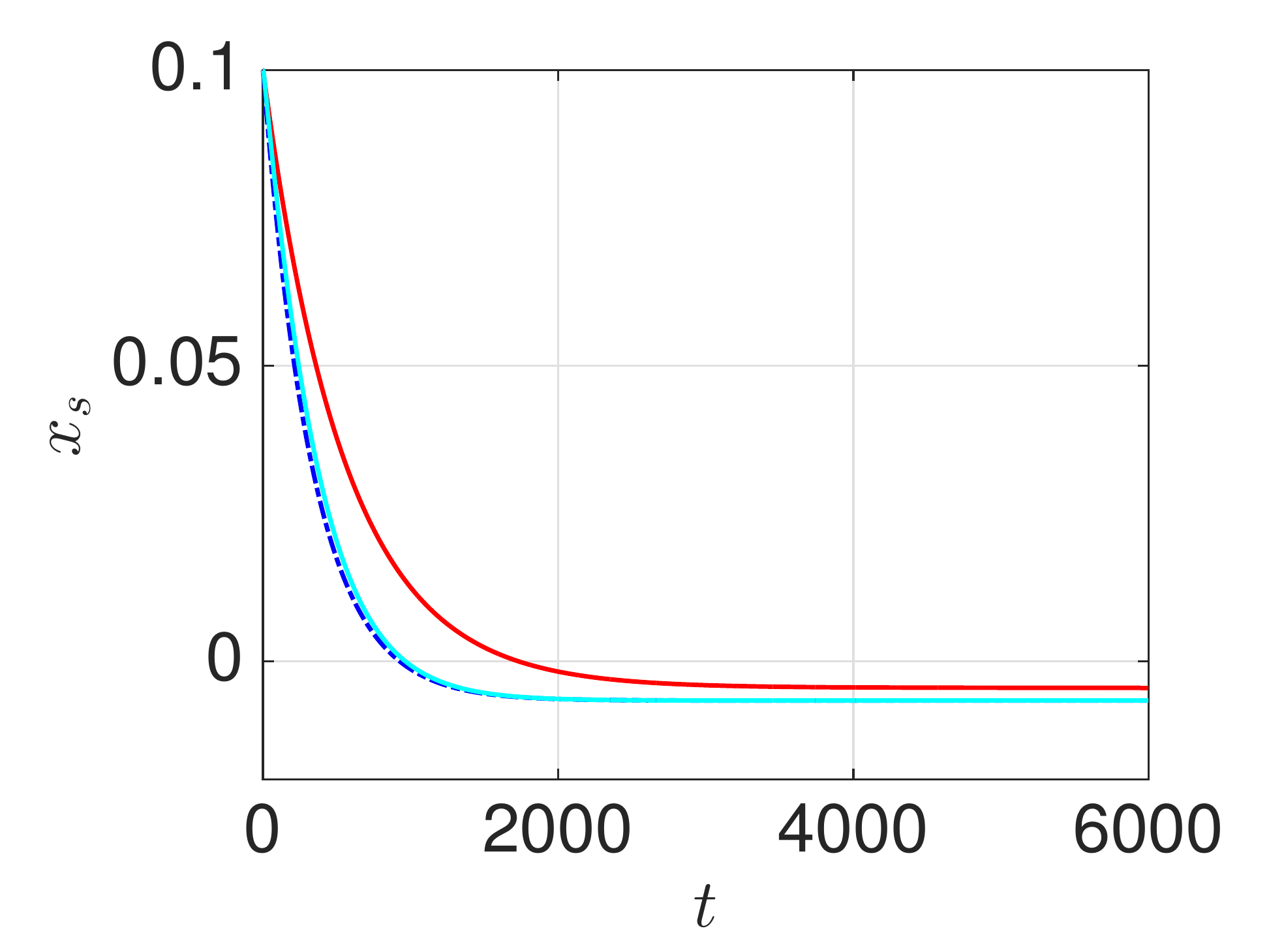}%
\caption{$\beta-\beta_c<0$}\label{simuxffix2}
\end{subfigure}
\quad
\begin{subfigure}[t]{0.47\columnwidth}
\centering
\includegraphics[width=\columnwidth,height=3.5cm]{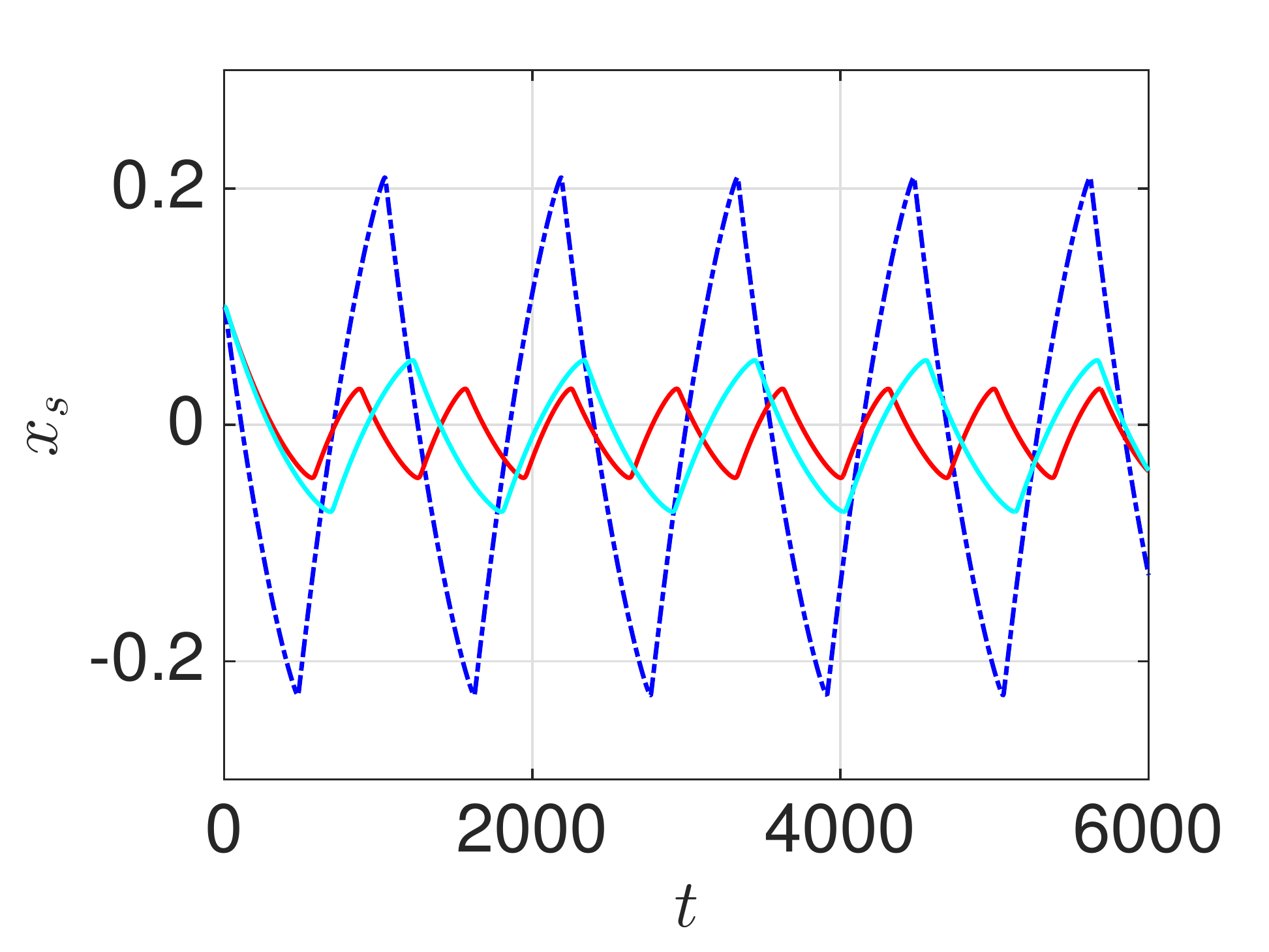}%
\caption{$0<\beta-\beta_c<1$}\label{simuxflc2}
\end{subfigure}%
\caption{States $x_f,x_s$ of system (\ref{EQ: sigmoidal hysteresis class}) for different nonlinearities. }\label{simupartdiff}
\end{figure}

We consider several sigmoid functions, namely $S_1:x\mapsto\tanh(x), S_2:x\mapsto e^{-e^{-x}}-e^{-1}, S_3:x\mapsto\frac{1}{1+e^{-x}}-\frac{1}{2}$. The corresponding critical values are $\beta_{c1}=1, \beta_{c2}=e$ and $\beta_{c3}=4$. We choose $\varepsilon=0.001, u=-0.01$ and $k=5$ in (\ref{est_xi}). Figure \ref{simupartdiff} illustrates the fact that, in spite of the different sigmoid nonlinearities, when $\beta-\beta_c<0$, the state $x_f$ converges to a constant value as time increases, as well as $x_s$, which is consistent with item \rm {1)} of Proposition \ref{LEM: stable osci dichotomy}. When $0<\beta-\beta_c<1$, both $x_f$ and $x_s$ converge to an oscillatory behavior, which is in agreement with item \rm {2)} of Proposition \ref{LEM: stable osci dichotomy}. 

\begin{figure}[H]
\centering
\begin{subfigure}[t]{0.47\columnwidth}
\centering
\includegraphics[width=\columnwidth,height=3.5cm]{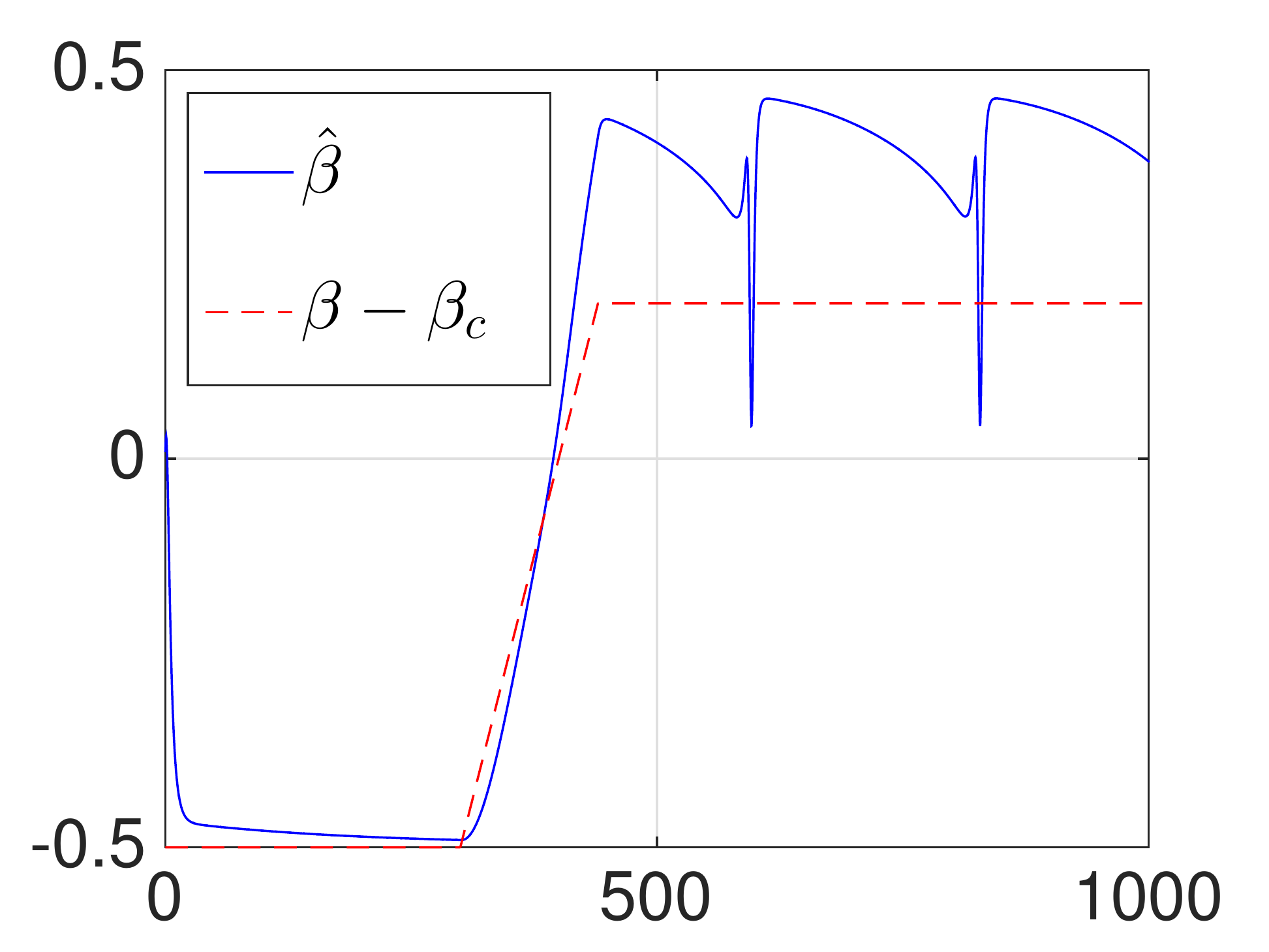}
\caption{without filter}\label{simupar}
\end{subfigure}
\quad
\begin{subfigure}[t]{0.47\columnwidth}
\centering
\includegraphics[width=\columnwidth,height=3.5cm]{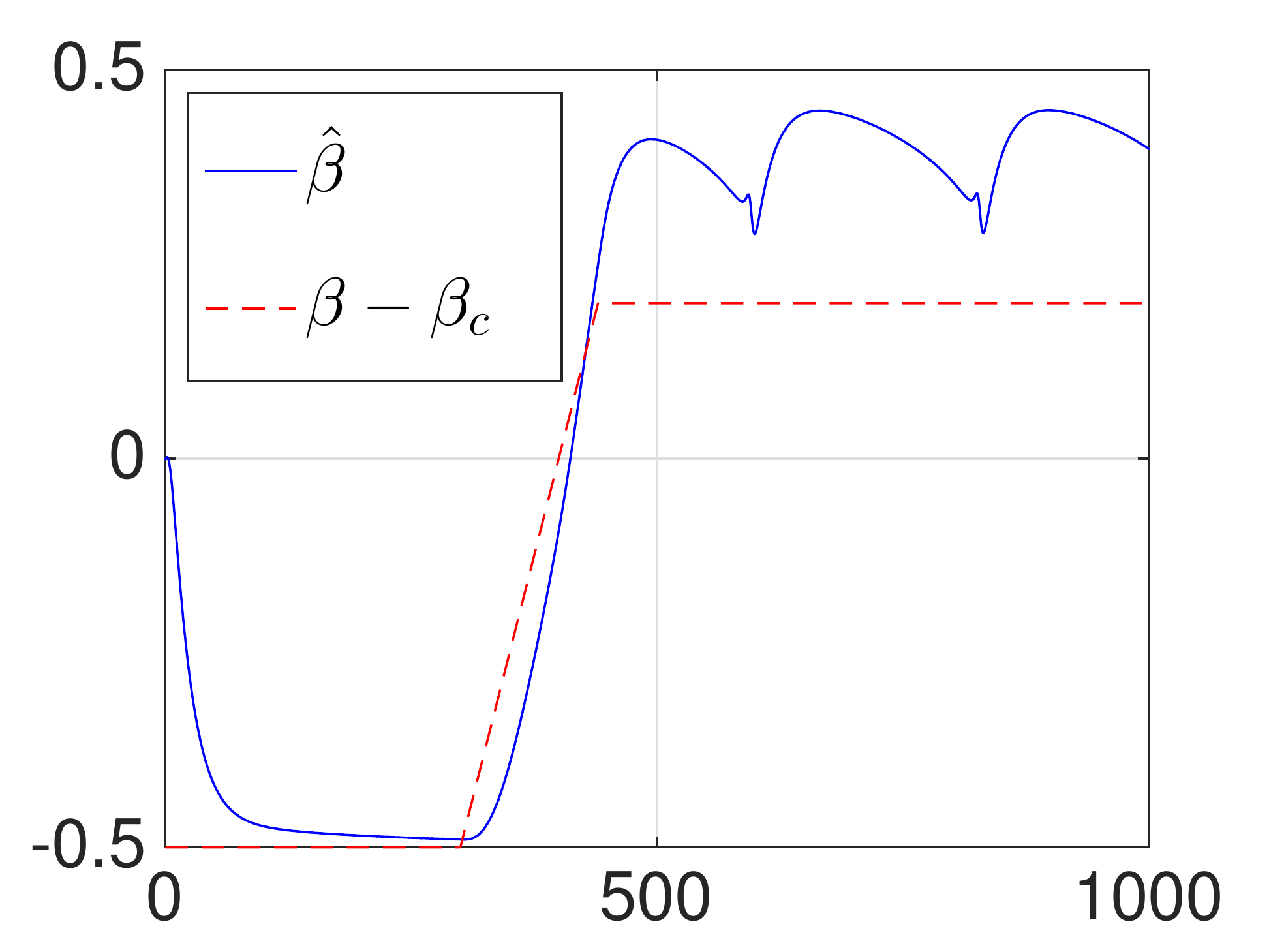}
\caption{with low pass filter}\label{simuparfilter}
\end{subfigure}
\caption{Variable $\hat\beta$ for two different values of $\beta$. }\label{simupart}
\end{figure}

To check that the proposed algorithm is able to detect a change of activity of the system, we use a varying signal for $\beta$ and we consider $S=S_1$ and the same data as above. In particular, $\beta(t)=0.5$ on $t\in[0,300]$,
it increases with slope $0.005$ on $t\in[300,425]$,
and $\beta(t)=1.2$ on $t\in[425,1000]$, which leads to a change of sign of $\beta-\beta_c$ at $t=400$, since $\beta=\beta_{c1}=1$ in this case. We observe in Figure \ref{simupar} that, when $\beta-\beta_c<0$, $\hat\beta$ converges to a constant value as time increases, corresponding to the resting activity.
When $0<\beta-\beta_c<1$, $\hat\beta$ tends to a periodic function, which is strictly positive. Hence, it indicates the oscillation activity. The spikes seen in Figure \ref{simupar} are due to the ``jump" of the solution to system (\ref{EQ: sigmoidal hysteresis class}) from one stable branch of the critical manifold to the other. This phenomenon is captured by item 2) of Theorem \ref{thm_betah}, as $\hat\beta$ is guaranteed to be close to $\hat \beta^*$ for all time except, periodically, over interval of length $|\delta T^\varepsilon|= \mathcal O(\varepsilon) T^\varepsilon$. These spikes can simply be removed by using a low pass filter on $\hat\beta$ as illustrated in Figure~\ref{simuparfilter}. Moreover, the value of $\hat\beta$ evolves from negative to positive when the change of sign of $\beta-\beta_c$ occurs. In other words, the qualitative detector is able to detect the current activity and whether a change is occuring. We have tested different values of $k$ in (\ref{est_xi}). The simulations indicate that the speed of convergence of $\hat\beta$ increases with $k$. However, this also leads to bigger spikes at ``jumps", which may provide wrong information during a very short interval, which may again be moderated by a low pass filter.

\begin{figure}[H]
\centering
\includegraphics[width=1\columnwidth,height=6cm]{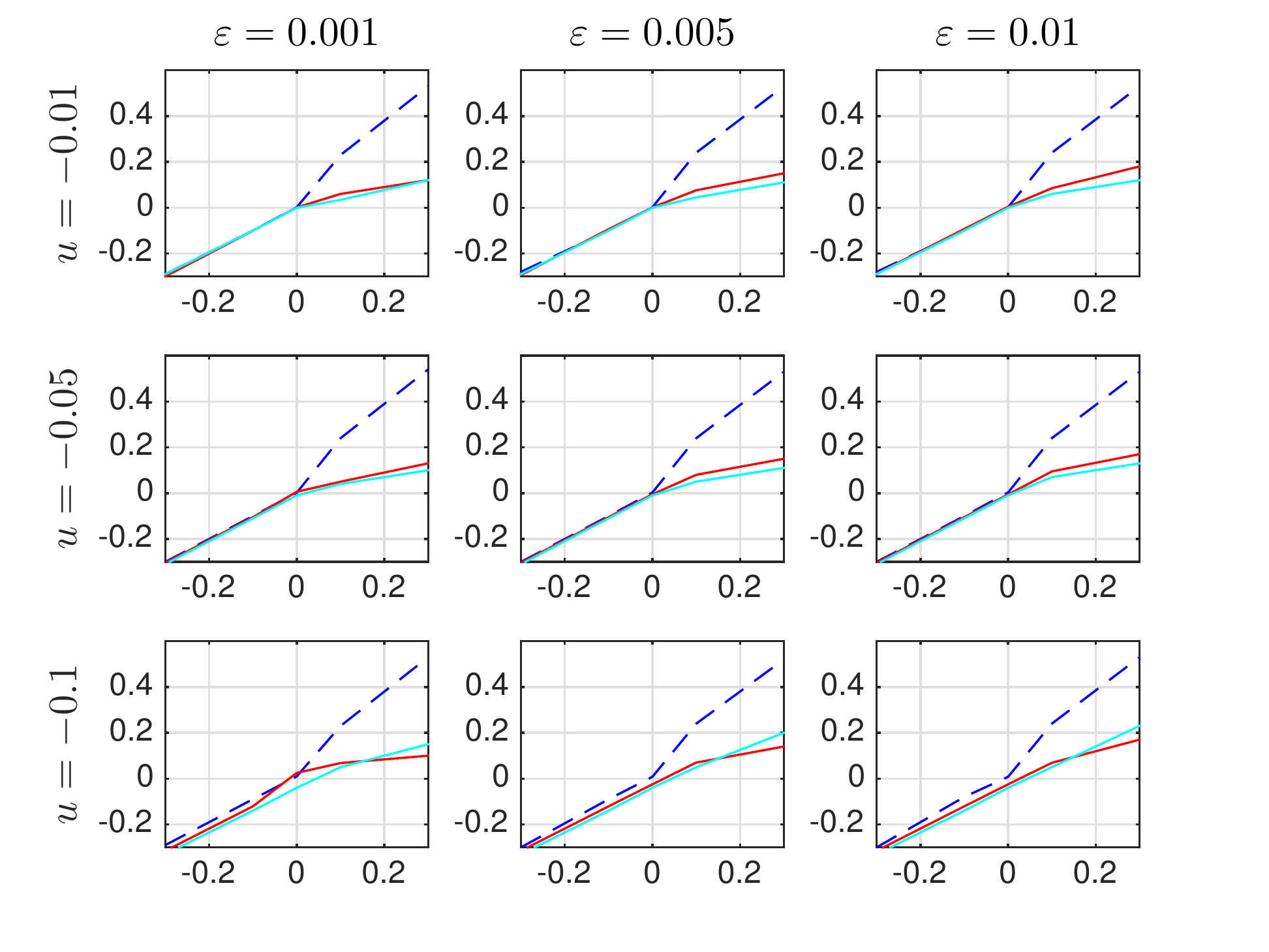}%
\caption{Evolution of $\hat\beta$ with $\beta-\beta_c$ for different sigmoid functions. Vertical axis is $\hat\beta$ and horizontal axis is $\beta-\beta_c$. Different $S$ functions: $S_1$ (blue dashed line), $S_2$ (red solid line), $S_3$ (light blue solid line).}
\label{simu1}
\end{figure}

The same tests have been done for $S_2$ and $S_3 $. We emphasize that, even though the nonlinearities are different, the detector remains the same as defined in (\ref{est_xi}). Figure \ref{simu1} shows the relationship between $\hat\beta$ and $\beta-\beta_c$ for different input $u$ and perturbation parameter $\varepsilon$. When $\beta-\beta_c<0$, the value of $\hat\beta$ provided in Figure \ref{simu1} is the constant value to which it converges as seen in simulations. When $0<\beta-\beta_c<1$, the value of $\hat\beta$ reported in the figure corresponds to the average value of $\hat\beta$. 
We observe that for $\beta-\beta_c<0$, $\hat\beta$ is negative and it increases as $\beta-\beta_c$ increases. Moreover, because the fixed point of system (\ref{EQ: sigmoidal hysteresis class}) is near zero for small $u$, $\hat\beta$ evolves almost linearly with respect to $\beta-\beta_c$.
 When $0<\beta-\beta_c<1$, $\hat\beta$ is positive. In addition, when $\beta-\beta_c=0$, $\hat\beta$ is around the origin. Hence $\hat\beta$ is able to detect the current activity type and whether a change is likely to occur. Figure \ref{simu1} illustrates the efficiency of the approach.

To further evaluate the robustness of the scheme, we have added small additive measurement noises given by $d_{out}(t)=0.005\sin(50t)$ for $x_f, x_s$ and $d_{in}=0.008$ for input $u$. It is found from Figure \ref{simurub} that the detector still provides good results in this case, which is in agreement with Remark \ref{rubst}.

\subsection{Hodgkin-Huxley model}\label{sec4.2}

We consider Hodgkin-Huxley model \cite{HH1952}
\begin{subequations}\label{HH model}
\ \begin{empheq}[left=\empheqlbrace]{align}
\label{hh1}
C\dot V&= \begin{aligned}[t]
&-g_Kn^4(V-V_K)\\
&-g_{Na}m^3h(V-V_{Na})\\
&-g_l(V-V_l)+I_{app},
\end{aligned}\\\label{hh2}
\tau_m(V)\dot m&=-m+m_\infty(V),\\\label{hh3}
\tau_n(V)\dot n&=-n+n_\infty(V),\\\label{hh4}
\tau_h(V)\dot h&=-h+h_\infty(V),
 \end{empheq}
\end{subequations}
where the variable $V$ is the membrane potential. The sodium fast activation is $m$ and its slow inactivation is $h$. The potassium slow activation is $n$. The input is the applied current $I_{app}$. The parameters $V_K, V_{Na}$ are the equilibrium potentials for the potassium and sodium ions and $V_l$ is the potential at which the leakage current is zero. The constant $C$ is the membrane capacity. The sodium, potassium and ionic conductances are denoted by $g_{Na}, g_K,$ and $g_l$.  The gating variable time constants are $\tau_m, \tau_n, \tau_h$ and their steady state characteristics are $m_\infty, n_\infty, h_\infty$, which are all monotone sigmoidal functions. We consider the same parameter values as in \cite{HH1952} unless otherwise specified.

As explained in Remark \ref{remHH}, Hodgkin-Huxley model can be reduced to two dimensions by assuming that $m\equiv m_\infty(V)$ and that $h\equiv l(n)$, for a suitable linear function $l$, which is fine close to the Hopf bifurcation point via center manifold reduction \cite{Guckenheimer1983}, see e.g., \cite{fitzhugh1961,rinzel1985}. In its two-dimensional reduced form, Hodgkin-Huxley model falls into the class of systems~\eqref{EQ: generic nonlinear}, with $V$ the fast variable and $n$ the slow variable. The fast nullcline can also be shown to be either monotone or cubic depending on $g_{Na}$ and, since $n_\infty$ is monotone, the slow variable nullcline can locally be approximated as being linear. We therefore expect the detector designed for the class of system (\ref{EQ: sigmoidal hysteresis class}) to work locally for the Hodgkin-Huxley model with $u=I_{app}$ and $\beta=g_{Na}$.

In the following, we apply the proposed algorithm to the full Hodgkin-Huxley model in (\ref{HH model}) to detect changes due to the  ruling parameter $g_{Na}$. We assume that the state variables $V, n$ are measurable. The input $I_{app}$ is a constant and is assumed to be known. From (\ref{est_xi}), the detector for the Hodgkin-Huxley model is
\begin{equation}
\dot{\hat g}_{Na}=-kV(-V^3+\hat g_{Na} V+I_{app}-m),
\end{equation}
where $k>0$.

\begin{figure}[H]
\centering
\includegraphics[width=0.8\columnwidth,height=5cm]{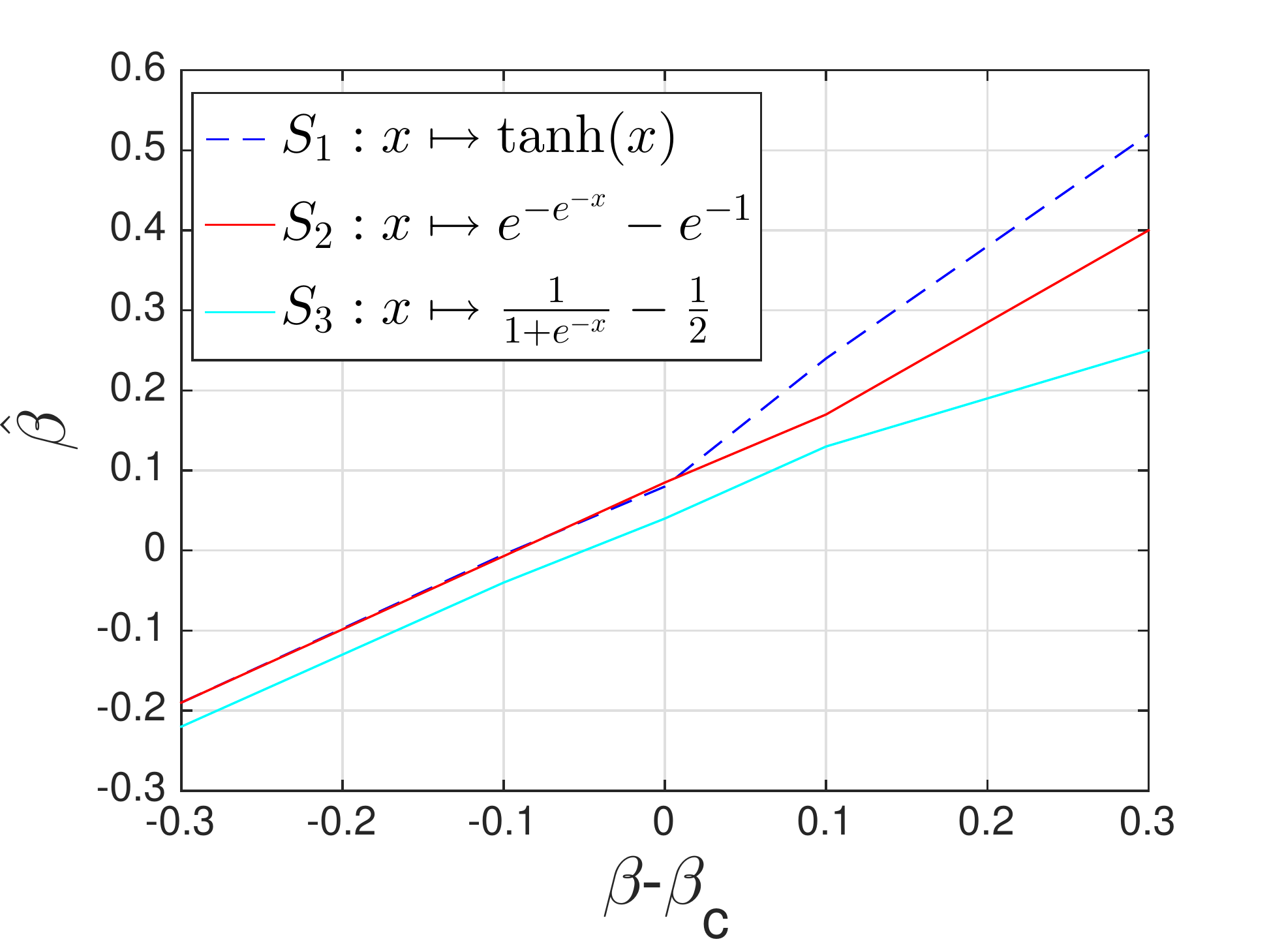}%
\caption{Evolution of $\hat\beta$ with $\beta-\beta_c$ when adding measurement noise.}
\label{simurub}
\end{figure}

\begin{figure}[H]
\centering
\includegraphics[width=0.8\columnwidth,height=5cm]{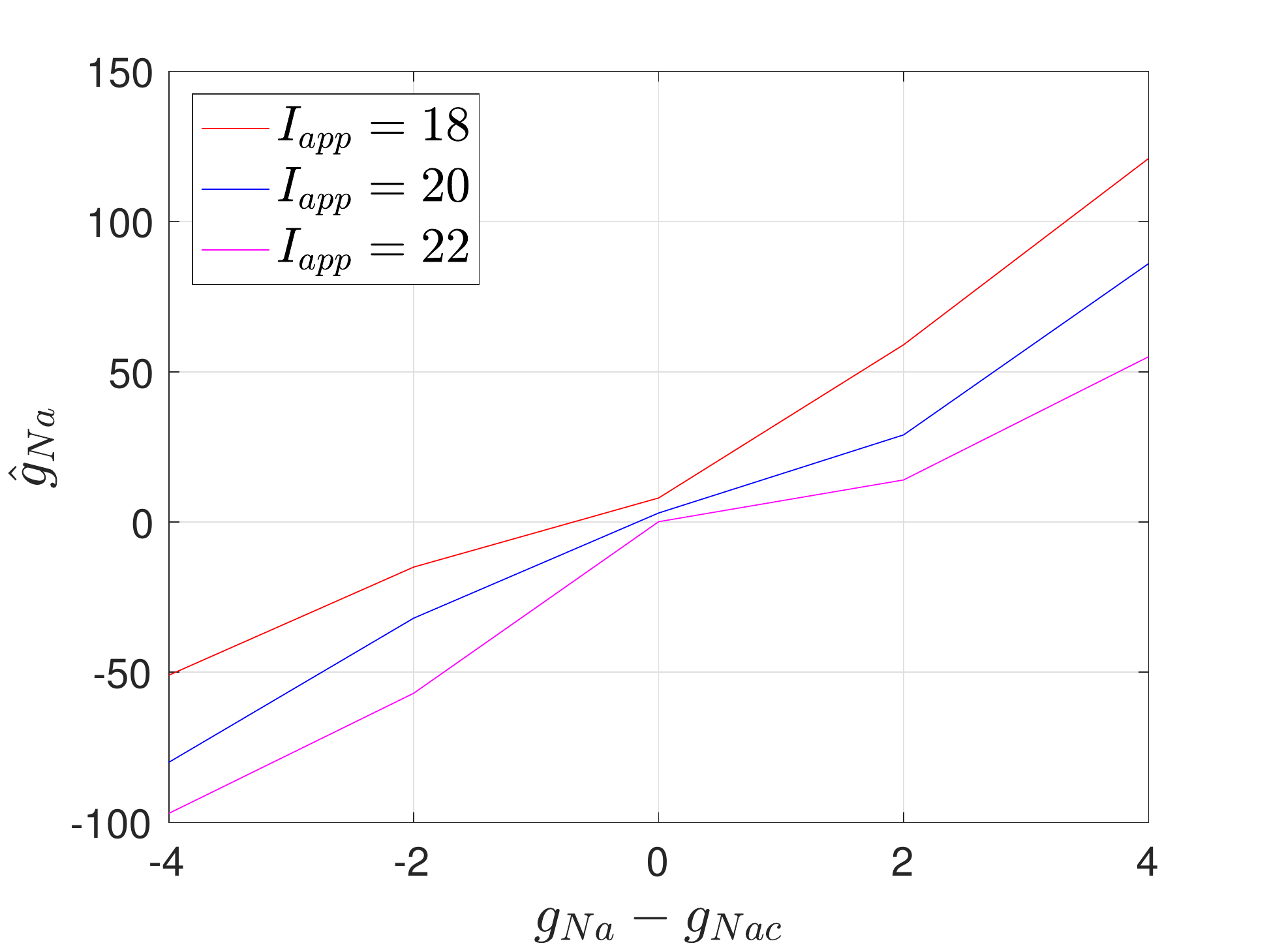}%
\caption{Evolution of the estimate $\hat g_{Na}$ with $g_{Na}-{g}_{Nac}$ for different input values $I_{app}$.}
\label{robustHH}
\end{figure}

In the simulations, the state variables $V, m$ are centered at zero by using a high-pass filter. We choose $k=4, I_{app}=20$. The values of $g_{Na}$ are selected from the set $\{23, 25,27,29,31\}$, which gives $g_{Na}-g_{Nac}\simeq\{-4,-2,0,2,4\}$ corresponding to an approximate critical value $g_{Nac}\simeq 27$. Figure~\ref{robustHH} shows the asymptotic value of $\hat g_{Na}$ as a function of $g_{Na}-g_{Nac}$, just like we did to generate Figures \ref{simu1}-\ref{simurub} in Section \ref{sec4.1}. The blue line corresponds to the nominal case $I_{app}=20$. Red and purple lines correspond to the perturbed cases $I_{app}=20\pm2$.
When $g_{Na}-g_{Nac}>0$, $\hat g_{Na}$ is selected as the average  value over a period of the asymptotic periodic function to which it converges. We observe that in the nominal case, $\hat g_{Na}$ is negative for $g_{Na}-g_{Nac}<0$ and it increases as $g_{Na}-g_{Nac}$ increases. The value of $\hat g_{Na}$ is positive when $g_{Na}-g_{Nac}>0$. Moreover, $\hat g_{Na}$ is around zero when $g_{Na}=g_{Nac}$.  Hence $\hat g_{Na}$ provides information about the actual activity of the model and whether the model is prone to change its activity. The perturbed cases also show that the detector performance is robust to small input uncertainties. These preliminary numerical results highlight the potentiality of the proposed approach beyond the academic example thoroughly analyzed in the present paper.

\section{Conclusion}\label{sec:5}
We have introduced the concept of {\it qualitative detection}, as the problem of informing on-line the qualitative dynamical behavior of a multiple-time scale nonlinear system, independently of large uncertainties on the system nonlinearities and without using any quantitative fitting of measured data. This is achieved by first extracting the system ruling parameter(s) and by subsequently designing a qualitative detector, which determines the system activity type and whether a qualitative change in the activity is close to occur.
We have presented this idea on a general class of nonlinear singularly perturbed systems. As a first application, we have focused on a class of two-dimensional nonlinear systems with two time-scales and a single nonlinearity either exhibiting resting or relaxation oscillation behaviors. We have illustrated the extension of the proposed detector to real-word application through numerical simulations. 
Future extensions will include generalizing the detector design problem to systems exhibiting more than two possible qualitatively distinct behaviors.

\section*{Appendix}\label{sec_appen}

{\bf Proof of Theorem \ref{property_quali_betahs}.} $1)$ Let $\beta\in \RR$ and $x_f\in S(\RR)\setminus \{0\}$. According to (\ref{ift_sin2}), $\hat\beta^*(x_f,\beta)=x^2_f-\frac{S^{-1}(x_f)}{x_f}+\beta$, which is well-defined and differentiable on $S(\RR)\setminus \{0\}\times \RR$. The only critical term for the smoothness of $\hat\beta^*$ with respect to $x_f$ is $\frac{S^{-1}(x_f)}{x_f}$ when $x_f=0$. In view of items a$)$ and c$)$ of Assumption \ref{sigmo_assum}, the function $S^{-1}$ is differentiable. Then the Taylor expansion of $S^{-1}(x_f)$ at the origin is, for any $x_f\in S(\RR)$
\begin{align}
\nonumber
S^{-1}(x_f)=&S^{-1}(0)+(S^{-1})^\prime(0)x_f\\\label{taylor_sinv}
&+\frac{(S^{-1})^{\prime\prime}(0)}{2}x^2_f+\mathcal O(x_f^3).
\end{align}
In view of item b$)$ of Assumption \ref{sigmo_assum} and (\ref{bc}), since $S(S^{-1}(0))=0$, it holds that 
$S^{-1}(0)=0$, $(S^{-1})^\prime(0)=\frac{1}{S^\prime (0)}=\beta_c>0$. Due to item d$)$ of Assumption \ref{sigmo_assum}, $(S^{-1})^{\prime\prime}(0)=-\frac{S^{\prime\prime}(0)}{(S^\prime (0))^2}=0$. We deduce from (\ref{taylor_sinv}), $S^{-1}(x_f)=\beta x_f+\mathcal O(x^3_f)$. Then the term $\frac{S^{-1}(x_f)}{x_f}$ becomes 
\begin{equation}
\label{taylor_sinv2}
\frac{S^{-1}(x_f)}{x_f}=\beta_c+\mathcal O(x^2_f),
\end{equation}
which is smooth with respect to $x_f$. This concludes the proof of the first item.\\
$2)$ Let $\beta\in \RR$, we first consider that $\beta<\beta_c$ and $x_f=x_s\in S^0$. By implicitly differentiating (\ref{sol_xs}) with respect to $\beta$, for $x_f\in S(\RR)$,
\begin{equation}
\label{der_beta_fix}
\frac{\partial x_f}{\partial \beta}=\frac{x_f }{-(\beta-1)+ (S^{-1})^\prime(x_f)}.
\end{equation}
Using (\ref{ift_sin2}) and (\ref{der_beta_fix}), we obtain
\begin{eqnarray}
\nonumber
\frac{\partial\hat\beta^*(x_f,\beta)}{\partial\beta}&=&\frac{\partial\hat\beta^*(x_f,\beta)}{\partial\beta}+\frac{\partial\hat\beta^*(x_f,\beta)}{\partial x_f}\frac{\partial x_f}{\partial\beta}\\\nonumber
&&\!\!\!\!\!\!\!\!\!\!\!\!\!\!\!\!\!\!\!\!\!\!\!\!\!\!\!\!=1+\frac{\left(2x_f-\frac{(S^{-1})^\prime(x_f)}{x_f}+\frac{S^{-1}(x_f)}{x^2_f}\right)x_f}{-(\beta-1)+ (S^{-1})^\prime(x_f)}\\
\label{termneg0}
&&\!\!\!\!\!\!\!\!\!\!\!\!\!\!\!\!\!\!\!\!\!\!\!\!\!\!\!\!=\frac{\left(1-\beta+2x^2_f+\frac{S^{-1}(x_f)}{x_f}\right)}{-(\beta-1)+ (S^{-1})^\prime(x_f)}. 
\end{eqnarray}
By item c$)$ of Assumption \ref{sigmo_assum}, $(S^{-1})^\prime(x_f)=\frac{1}{S^\prime(S^{-1}(x_f))}\geqslant \frac{1}{S^\prime(0)}=\beta_c$. Since $\beta-\beta_c<0$, we have $-(\beta-1)+(S^{-1})^\prime(x_f)\geqslant \beta_c-\beta+1>0$. Then the term $\frac{1}{-(\beta-1)+ (S^{-1})^\prime(x_f)}$ in the right-hand side of (\ref{termneg0}) is strictly positive.
According to the mean value theorem, it holds $\frac{S^{-1}(x_f)-S^{-1}(0)}{x_f-0}=\frac{S^{-1}(x_f)}{x_f}=(S^{-1})^\prime(c)$ for some $c$ such that $c\in(x_f,0)$ if $x_f<0$ and $c\in(0,x_f)$ if $x_f>0$. Using again item c$)$ of Assumption \ref{sigmo_assum}, it holds $(S^{-1})^\prime(c)=\frac{1}{S^\prime(c)} \geqslant \frac{1}{S^\prime(0)}= \beta_c$. Hence, we deduce from (\ref{termneg0}) with the fact that $\beta_c-\beta>0$, $\frac{\partial\hat\beta^*(x_f,\beta)}{\partial\beta}\geqslant\frac{(1+\beta_c-\beta+2x^2_f)}{-(\beta-1)+ (S^{-1})^\prime(x_f)}>0$.
\\
3$)$ Let $\beta\in(0,\beta_c+1)$, for any fixed $x_s\in \RR$ on $S^0$ and $u\in \RR$, by implicitly differentiating (\ref{sol_xs}) with respect to $\beta$, we have $\frac{\partial x_f}{\partial \beta}=\frac{x_f}{-\beta +(S^{-1})^\prime(x_f)}$,
which is well-defined as long as $x_f\neq \{x^-_{fold},x^+_{fold}\}$.
\\
We deduce from (\ref{ift_sin2}) and $\frac{\partial x_f}{\partial \beta}=\frac{x_f}{-\beta +(S^{-1})^\prime(x_f)}$,
\begin{eqnarray}
\nonumber
\frac{\partial\hat\beta^*(x_f,\beta)}{\partial\beta}&=&\frac{\partial\hat\beta^*(x_f,\beta)}{\partial\beta}+\frac{\partial\hat\beta^*(x_f,\beta)}{\partial x_f}\frac{\partial x_f}{\partial\beta}\\\label{dhbetastar}
&=&1+\frac{\left(2x^2_f-(S^{-1})^\prime(x_f)+\frac{S^{-1}(x_f)}{x_f}\right)}{-\beta +(S^{-1})^\prime(x_f)}.
\end{eqnarray}
We use again $\frac{S^{-1}(x_f)}{x_f}=(S^{-1})^\prime(c)$ for some $c$ such that $c\in(x_f,0)$ if $x_f<0$ and $c\in(0,x_f)$ if $x_f>0$, and write $(S^{-1})^\prime(c)=(S^{-1})^\prime(x_f)+(S^{-1})^\prime(c)-(S^{-1})^\prime(x_f)$. Substituting this expression into \eqref{dhbetastar}, we obtain
\begin{equation}\label{value_case2}
\frac{\partial\hat\beta^*(x_f,\beta)}{\partial\beta}= 1+\frac{(S^{-1})^\prime(c)-(S^{-1})^\prime(x_f)+2x_f^2}{-\beta+(S^{-1})^\prime (x_f)},
\end{equation}
since $(S^{-1})^\prime(c)=\beta_c+\mathcal O(c^2)$, $(S^{-1})^\prime(x_f)=\beta_c+\mathcal O(x^2_f)$ and $|c|<|x_f|$, it follows from (\ref{value_case2}) that $\frac{\partial\hat\beta^*(x_f,\beta)}{\partial\beta}=1+\mathcal O(x_f^2)$. This proves the first part of item 3$)$.

We now prove the second part of item 3$)$. Let $x_{fold}\in\{x^+_{fold}, x^-_{fold}\}$. By definition of $x^+_{fold}, x^-_{fold}$ in Lemma \ref{property_quali_betahs}, it holds $-1+\beta S^\prime(S^{-1}(x_{fold}))=0,$
which is equivalent to
\begin{equation}
\label{eq_fold_2}
\beta=(S^{-1})^\prime(x_{fold}),
\end{equation}
as $(S^{-1})^\prime(x_{fold})=\frac{1}{S^\prime(S^{-1}(x_{fold}))}$.

We deduce from (\ref{dhbetastar}) by using (\ref{eq_fold_2}), 
\begin{eqnarray}
\nonumber
\frac{\partial\hat\beta^*(x_f,\beta)}{\partial\beta}&=&1+\frac{\left(2x^2_f-(S^{-1})^\prime(x_f)+\frac{S^{-1}(x_f)}{x_f}\right)}{-\beta +(S^{-1})^\prime(x_f)}\\\nonumber
&=&\frac{\left(-\beta+2x^2_f+\frac{S^{-1}(x_f)}{x_f}\right)}{-\beta+(S^{-1})^\prime (x_f)}\\\label{termpos}
&=&\frac{\left(-(S^{-1})^\prime(x_{fold})+2x^2_f+\frac{S^{-1}(x_f)}{x_f}\right)}{-\beta+(S^{-1})^\prime (x_f)}.
\end{eqnarray}
By virtue of items c$)$ and d$)$ of Assumption \ref{sigmo_assum}, the following holds, for $x_f\in(-\infty,x^-_{fold})\cup (x^+_{fold},+\infty)$
\begin{equation}
\label{sta_branch_comp}
-1+\beta S^\prime(S^{-1}(x_f))<0.
\end{equation}
Dividing both sides of the above inequality by $S^\prime(S^{-1}(x_{f}))$, which is strictly positive, and using $(S^{-1})^\prime(x_{f})=\frac{1}{S^\prime(S^{-1}(x_{f}))}$, we obtain
\begin{align}
\label{term_311}
-\beta+(S^{-1})^\prime(x_f)&>0.
\end{align}
This implies that the denominator of (\ref{termpos}) is strictly positive. We denote $g(x_f)=\frac{(S^{-1}) (x_f)}{x_f}$. The derivative of $g(x_f)$ is given by
\begin{align}
\label{der_phi}
g^\prime(x_f)
&=\frac{x_f(S^{-1})^\prime (x_f)-S^{-1} (x_f)}{x^2_f}.
\end{align}
We show below that the growth of $g(x_f)$ increases with $\vert x_f\vert$. According to item c$)$ of Assumption \ref{sigmo_assum}, it holds $(S^{-1})^\prime(x_f) =\frac{1}{S^\prime(S^{-1}(x_f))}>0$. Hence, the function $S^{-1}$ is strictly increasing. Moreover, item b$)$ of Assumption \ref{sigmo_assum} implies $S^{-1}(0)=0$. We first consider the case where $x_f>0$. It holds that $S^{-1}(x_f)>0$ for all $x_f>0$. Using $ \frac{(S^{-1}) (x_f)}{x_f}=(S^{-1})^\prime (c)$ as previously done in this proof, for some $c$ such that $c\in(0,x_f)$ and in view of items c$)$ and d$)$ of Assumption \ref{sigmo_assum}, $S^\prime(S^{-1}(x_f))<S^\prime(S^{-1}(c))$ for $x_f>c>0$, and thus $(S^{-1})^\prime (x_f)=\frac{1}{S^\prime(S^{-1}(x_f))}>\frac{1}{S^\prime(S^{-1}(c))}=(S^{-1})^\prime (c)$. It holds that $\frac{(S^{-1}) (x_f)}{x_f}=(S^{-1})^\prime (c)<(S^{-1})^\prime (x_f)$. Then, we deduce
\begin{align}
\label{der_phi_num1}
&x_f(S^{-1})^\prime (x_f)-(S^{-1}) (x_f)>0.
\end{align}
In view of (\ref{der_phi}) and (\ref{der_phi_num1}), we have $g^\prime(x_f)>0$ for all $x_f>0$, which means that $\frac{S^{-1}(x_f)}{x_f}$ increases with $x_f$ when $x_f>0$. Similar arguments show that 
 $\frac{S^{-1}(x_f)}{x_f}$ increases with $ \vert x_f\vert$ when $x_f<0$. Therefore, the term $\frac{S^{-1}(x_f)}{x_f}$, as well as $2x^2_f$, increases when $x_f\in(-\infty,x^-_{fold})\cup (x^+_{fold},+\infty)$. Under condition (\ref{cond_suf}), $-(S^{-1})^\prime(x_{fold})+2x^2_f+\frac{S^{-1}(x_f)}{x_f}$ is strictly positive for $x_f\in(-\infty,x^-_{fold})\cup (x^+_{fold},+\infty)$. Recall that $\frac{1}{-\beta+(S^{-1})^\prime (x_f)}$ in (\ref{termpos}) is strictly positive, we obtain the desired result.

4$)$ Combining (\ref{ift_sin2}) and (\ref{taylor_sinv2}), we obtain for any $x_f\in S(\RR)$ and $\beta\in\RR$, $\hat\beta^*(x_f,\beta)= \beta-\beta_c+\mathcal O (x^2_f)$.
\\
This concludes the proof of the last item. $\hfill\blacksquare$

{\bf Proof of Proposition \ref{proposition2}.} 
Two steps are used to prove this proposition. We first show a boundedness property for system (\ref{est_xi}). Then we prove that (\ref{diss_eq}) holds. Let $\Delta,M,T,\mu>0$ and $x_f,x_s,u$ be such that conditions $(i)$ and $(ii)$ of Proposition \ref{proposition2} hold.
We rewrite (\ref{est_xi}) as follows
\begin{equation}
\label{est_xi_tran}
\dot {\hat \beta}=-kx^2_f\hat\beta+\biggl(kx^4_f+kx_f(x_s-u)\biggr).
\end{equation}
System (\ref{est_xi_tran}) can be interpreted as a linear time-varying system with input $kx^4_f+kx_f(x_s-u)$. We first consider the following nominal system 
\begin{equation}
\label{est_xi_tran_nom}
\dot {\hat \beta}=-kx^2_f\hat\beta.
\end{equation}
Solutions to system (\ref{est_xi_tran_nom}) are given for $t\geqslant 0$ by $\hat\beta(t)=\Phi(t,0)\hat\beta(0)$,
where $\Phi(t,s):=\exp\left(\int_s^t -k x^2_f(\tau)\; d\tau\right)$ and $\vert \Phi(t,0)\vert \leqslant 1$ for all $t\geqslant 0$.
Hence the origin is stable for system (\ref{est_xi_tran_nom}).\\
Let $\vert \hat\beta(0)\vert\leqslant \Delta$ and consider the following Lyapunov function candidate for system (\ref{est_xi_tran}) 
\begin{equation}
\label{lya_bound}
V(t,\hat\beta):=p(t)\hat\beta^2,
\end{equation}
where $p(t)\!\!:=\!\!\int_t^\infty \Phi^2(s,t)\; \!\!ds$. We show below that $V(t,\hat\beta)$ is positive definite and radially unbounded with respect to $\hat\beta$, uniformly in $t$. For $0\leqslant t\leqslant s$, we write
\begin{equation}
\label{phi_gen}
p(t)=\int_t^\infty \Phi^2(s,t)\; ds=\int_t^\infty \exp\left( \int_t^s -2k x^2_f(\tau)\; d\tau\right) ds.
\end{equation}
Since $x_f$ is $(T,\mu)$-PE by assumption, for $t\leqslant s< t+T$, $\mu\leqslant \int_t^{t+T} x^2_f(\tau) d\tau=\int_t^s x^2_f(\tau) d\tau+\int_s^{t+T} x^2_f(\tau) d\tau$. Hence,
\begin{eqnarray}
\label{cond1cha}
\int_t^s x^2_f(\tau) d\tau &\geqslant& \mu-\int_s^{t+T} x^2_f(\tau) d\tau.
\end{eqnarray}
Under assumption $\max(\Vert x_{f}\Vert_\infty,\Vert x_{s}\Vert_\infty,\Vert u\Vert_\infty)\leqslant M$, we deduce from (\ref{phi_gen}) and (\ref{cond1cha}) 
\begin{eqnarray}
\nonumber
p(t)&=&\int_t^\infty \exp\left( \int_t^s -2k x^2_f(\tau)\; d\tau\right) ds\\\nonumber
&\leqslant& \int_t^\infty \exp\left( -2k\mu+2k\int^{t+T}_s x^2_f(\tau)\; d\tau\right) ds\\\nonumber
&\leqslant& \exp(-2k\mu)\int_t^\infty \exp\left( 2kM^2(t+T-s)\right) ds\\\label{phi_up}
&\leqslant&\frac{\exp(2kM^2T)}{2kM^2}.
\end{eqnarray}
We also have from (\ref{phi_gen})
\begin{align}
\nonumber
p(t)&=\int_t^\infty \exp\left( \int_t^s -2k x^2_f(\tau)\; d\tau\right) ds\\\label{phi_down}
&\geqslant \int_t^\infty \exp\left(  -2k M^2(s-t)\right) ds \geqslant \frac{1}{2kM^2}.
\end{align}
Hence, $V(t,\hat\beta)$ is lower and upper bounded as
\begin{equation}
\label{p_bound}
\frac{1}{2kM^2}\vert \hat\beta\vert^2\leqslant V(t,\hat\beta) \leqslant \frac{\exp(2kM^2T)}{2kM^2}\vert \hat\beta\vert^2.
\end{equation}
The time-derivative of $V(t,\hat\beta)$ along the solution to (\ref{est_xi_tran}) is, for any $t\geqslant 0$,
\begin{eqnarray}
\nonumber
&&\frac{\partial V}{\partial t}+\frac{\partial V}{\partial \hat\beta} \hat f(\hat\beta, x_f, u-x_s)\\\nonumber
&=&\dot p(t)\hat{\beta}^2+2p(t)\hat\beta\left(-kx^2_f\hat\beta+\left(kx^4_f+kx_f(x_s-u)\right)\right)\\\nonumber
&=&\left(\dot p(t)-2kx^2_fp(t)\right)\hat\beta^2+2\hat\beta p(t)\left(kx^4_f+kx_f(x_s-u)\right).\\\label{dv_1}
\end{eqnarray}
Since $\Phi(t,t)=1$ and $\frac{d\Phi(s,t)}{dt}=\Phi(s,t)kx^2_f$, we deduce after some calculations that
\begin{equation}
\label{dP_fin}
\dot p(t)-2kx^2_fp(t)=-1.
\end{equation}
Using (\ref{phi_up}) and (\ref{dP_fin}), $\frac{\partial V}{\partial t}+\frac{\partial V}{\partial \hat\beta} \hat f(\hat\beta, x_f, u-x_s)\leqslant- \vert \hat\beta\vert^2+2\vert \hat\beta\vert \frac{\exp(2kM^2T)\tilde C(M)}{2kM^2}$.
According to \cite[Theorem 4.18]{Nlsys3rd}, we conclude the boundedness property of (\ref{est_xi}) from (\ref{p_bound}) and $\frac{\partial V}{\partial t}+\frac{\partial V}{\partial \hat\beta} \hat f(\hat\beta, x_f, u-x_s)\leqslant- \vert \hat\beta\vert^2+2\vert \hat\beta\vert \frac{\exp(2kM^2T)\tilde C(M)}{2kM^2}$. That is, there exists $\bar\eta(\Delta, M,T)>0$ such that solution $\hat\beta$ to system (\ref{est_xi_tran}) holds $\vert \hat\beta(t) \vert \leqslant\bar \eta(\Delta,M,T)$, for all $t\geqslant 0$.
\\
To prove (\ref{diss_eq}), let $x_f,x_s,u$ be such that conditions $(i)$ and $(ii)$ of Proposition \ref{proposition2} holds, we define for $t\geqslant 0$, 
\begin{equation}
\label{lya_auxi}
W(t,\hat\beta_1,\hat\beta_2):=p_1(t)(\hat\beta_1-\hat\beta_2)^2,
\end{equation}
where $p_1(t):=\int_t^\infty \Phi_1^2(s,t)\; ds$, and $\Phi_1(t,0):= \exp\left(\int_0^{t} -k x^2_{f1}(\tau) d\tau\right)$. For $\vert \hat\beta_1(0)\vert, \vert \hat\beta_2(0)\vert\leqslant \Delta$, similar computations in (\ref{phi_up})-(\ref{phi_down}) shows, for any $t\geqslant 0$
\begin{eqnarray}
\nonumber
\frac{1}{2kM^2} \vert \hat\beta_1-\hat\beta_2\vert^2&\leqslant& W(t,\hat\beta_1,\hat\beta_2)\\\label{bound_lyapunow}
&\leqslant&\frac{\exp(2kM^2T)}{2kM^2}\vert \hat\beta_1-\hat\beta_2\vert^2.
\end{eqnarray}
Following the similar computations for (\ref{dP_fin}), we obtain $\dot p_1(t)-2kx^2_{f1}p_1(t)=-1$.
The following holds along to the solutions to (\ref{est_xi_tran}) for any $t\geqslant 0$, where the time arguments are omitted,
\begin{eqnarray}
\nonumber
&&\!\!\!\!\!\!\!\!\frac{\partial W}{\partial t}+\frac{\partial W}{\partial \hat\beta_1} \hat f(\hat\beta_1, x_{f1}, u_1-x_{s1})+\frac{\partial W}{\partial \hat\beta_2} \hat f(\hat\beta_2, x_{f2}, u_2-x_{s2})\\\nonumber
&\!\!\!\!=&\!\!\!\dot p_1(t)(\hat \beta_1-\hat\beta_2)^2+2p_1(t)(\hat\beta_1-\hat\beta_2)\biggl(-kx^2_{f1}\hat\beta_1+\left(kx_{f1}^4\right.\\\nonumber
&&\!\!\!\!\!\!\!\!\left.+kx_{f1}(x_{s1}-u_1)\right)+kx^2_{f2}\hat\beta_2-\left(kx_{f2}^4+kx_{f2}(x_{s2}-u_2)\right)\biggr)\\\nonumber
&=& \left(\dot p_1(t)-2kx^2_{f1}p_1(t)\right)(\hat \beta_1-\hat\beta_2)^2+2 kp_1(t)(\hat\beta_1-\hat\beta_2)\\\nonumber
&&\!\!\!\!\times\biggl(-\hat\beta_2 (x^2_{f1}-x^2_{f2})+(x^4_{f1}-x^4_{f2})+x_{s2}(x_{f1}-x_{f2})\\\label{dw1_1est}
&&\!\!\!\!+x_{f1}(x_{s1}-x_{s2})-x_{f1}(u_1-u_2)-u_2(x_{f1}-x_{f2})\biggr).
\end{eqnarray}
Under condition $(ii)$, $\dot p_1(t)-2kx^2_{f1}p_1(t)=-1$, and the fact that $x \mapsto x^4$ and $x\mapsto x^2$ are locally Lipschitz, we derive that there exists $\rho(\Delta,M,T)>0$ such that $\frac{\partial W}{\partial t}+\frac{\partial W}{\partial \hat\beta_1} \hat f(\hat\beta_1, x_{f1}, u_1-x_{s1})+\frac{\partial W}{\partial \hat\beta_2} \hat f(\hat\beta_2, x_{f2}, u_2-x_{s2})\leqslant-\vert \hat \beta_1-\hat\beta_2\vert^2+ \rho(\Delta,M,T)\vert \hat\beta_1-\hat\beta_2\vert (\vert x_{f1}-x_{f2}\vert+\vert x_{s1}-x_{s2}\vert+\vert u_1-u_2\vert )$,
we then obtain the desired result by following the similar analysis as in the first step of this proof. $\hfill\blacksquare$

{\bf Proof of Lemma \ref{lem_PE}.} We first prove item \rm{1)}. Let $\Delta>\delta>0$. For any $\beta-\beta_c<0$, any constant input $\delta<\vert u\vert<\Delta$, according to item \rm{1} of Proposition \ref{LEM: stable osci dichotomy}, there exists $\bar \varepsilon$ such that for any $\varepsilon\in(0,\bar\varepsilon]$ and any $\vert x_f(0)\vert, \vert x_s(0)\vert\leqslant \Delta$, it holds $\vert (x_f(t)-x_f^*),(x_s(t)-x_s^*)\vert\leqslant c(\Delta,\delta)e^{-\sigma(\Delta,\delta)t}\vert (x_f(0)-x^*_f),(x_s(0)-x^*_s)\vert$, for $t\geqslant 0$, where the fixed point $(x^*_f,x^*_s)$ is different from $(0,0)$ (because $u\neq 0$) and $c(\Delta,\delta), \sigma(\Delta,\delta)>0$. 
Hence, condition $(ii)$ of item 1) of Lemma \ref{lem_PE} holds. We also deduce that there exists $T_1(\Delta,\delta)>0$ such that $\vert x_f(t)-x_f^*\vert\leqslant \frac{\vert x_f^*\vert}{2}$ for any $t\geqslant T_1(\Delta,\delta)$. It holds that for $t\geqslant T_1(\Delta,\delta)$, $\vert x_f(t)\vert=\vert x_f(t)-x^*_f+x^*_f\vert\geqslant \vert x^*_f\vert -\vert x_f(t)-x^*_f\vert\geqslant \vert x^*_f\vert -\frac{\vert x^*_f\vert}{2}= \frac{\vert x^*_f\vert}{2}$.
Let $T(\Delta,\delta)>T_1(\Delta,\delta)$, we deduce, for any $t\geqslant 0$
\begin{align}
\nonumber
&\int_t^{t+T(\Delta,\delta)}\!\!\!\! x^2_f(\tau) d\tau=\!\!\int_t^{t+T_1(\Delta,\delta)}\!\!\!\! \!\!x^2_f(\tau) d\tau+\!\!\int_{t+T_1(\Delta,\delta)}^{t+T(\Delta,\delta)} \!\!\!\!\!\!x^2_f(\tau) d\tau\\\nonumber
&\geqslant \int_{t+T_1(\Delta,\delta)}^{t+T(\Delta,\delta)} \vert x_f(\tau)\vert^2 d\tau \geqslant \left(\frac{\vert x^*_f\vert}{2}\right)^2\!\!\!\!(T(\Delta,\delta)-T_1(\Delta,\delta)).
\end{align}
Thus, condition $(i)$ of item 1) of Lemma \ref{lem_PE} holds with $\mu(\Delta,\delta)=\left(\frac{\vert x^*_f\vert}{2}\right)^2\times(T(\Delta,\delta)-T_1(\Delta,\delta))$. 

We next prove item 2). Let $\Delta>\delta>0$, for any $0<\beta-\beta_c<1$ and for any input $u\in(-\bar u,\bar u)$, where $\bar u\in(0,\Delta]$, according to item \rm{2)} of Proposition \ref{LEM: stable osci dichotomy}, there exists $\bar \varepsilon$ such that for all $\varepsilon\in(0,\bar\varepsilon]$, and $\vert x_f(0)\vert, \vert x_s(0)\vert\leqslant \Delta$ and $\vert x_f(0)-p^*_f\vert, \vert x_s(0)-p^*_s\vert\geqslant \delta$, all trajectories of system (\ref{EQ: sigmoidal hysteresis class}) converge to an asymptotic periodic orbit $P^\varepsilon$ with period $T^\varepsilon$ and it is locally exponentially stable. 
Thus solutions $x_f,x_s$ to system (\ref{EQ: sigmoidal hysteresis class}) are bounded. Therefore, condition $(ii)$ of item 2) of Lemma \ref{lem_PE} holds. Moreover, let $\vert x_f^{lc}\vert$ be the absolute value of $x_f$-component of $P^\varepsilon$, which is $T^\varepsilon$-periodic. The average value of $\vert x^{lc}_f\vert$ over a period is denoted by 
$\bar x_f=\frac{1}{T^\varepsilon}\int_t^{t+T^{\varepsilon}} \vert x^{lc}_f(\tau)\vert d\tau$ for all $t\geqslant 0$ and it is strictly positive.
Due to the globally attractivity of $P^\varepsilon$, for any $\eta>0$ there exists $T_1(\Delta,\delta,\eta)>0$ such that $\vert x_f(t)-\bar x_f\vert\leqslant \eta$ for $t\geqslant T_1(\Delta,\delta,\eta)$. 
Let $T(\Delta,\delta)=T_1(\Delta,\delta,\eta)+T^\varepsilon$, the following holds 
\begin{eqnarray}
\nonumber
&&\int_t^{t+T(\Delta,\delta)}\!\!\!\!\!\! x^2_f(\tau) d\tau\!\!=\!\!\!\!\int_t^{t+T_1(\Delta,\delta,\eta)} \!\!\!\!\!\!\!\!x^2_f(\tau) d\tau+\!\! \int_{t+T_1(\Delta,\delta,\eta)}^{t+T(\Delta,\delta)} \!\!\!\!\!\!x^2_f(\tau) d\tau\\\nonumber
&&\geqslant  \int_{t+T_1(\Delta,\delta)}^{t+T(\Delta,\delta)} \bar x^{2}_f d\tau+\int_{t+T_1(\Delta,\delta)}^{t+T(\Delta,\delta)} (x^2_f(\tau)-\bar x^{2}_f) d\tau\\\label{cal_inter}
&&\geqslant   \bar x^{2}_f T^\varepsilon-\int_{t+T_1(\Delta,\delta)}^{t+T(\Delta,\delta)} \vert x^2_f(\tau)-\bar x^{2}_f\vert d\tau.
\end{eqnarray}
Since $x\mapsto x^2$ is locally Lipschitz, it holds $ \vert x^2_f(t)-\bar x^{2}_f\vert\leqslant \psi(\Delta,\delta)\vert x_f(t)-\bar x_f\vert$.
 We deduce from (\ref{cal_inter})
 \begin{eqnarray}
 \nonumber
 \int_t^{t+T(\Delta,\delta)} x^2_f(\tau) d\tau&\geqslant&  \bar x^{2}_f T^\varepsilon-\int_{t+T_1(\Delta,\delta)}^{t+T(\Delta,\delta)} \vert x^2_f(\tau)-\bar x^{2}_f\vert d\tau\\\label{cal_inter2}
 &&\geqslant  \bar x^{2}_f T^\varepsilon-\psi(\Delta,\delta)\eta T^\varepsilon.
 \end{eqnarray}
 Let choose $T_1(\Delta,\delta,\eta)$ sufficiently large such that $\eta=\frac{\bar x_f}{2\psi(\Delta,\delta)}$. Then, the above inequality follows
  \begin{align}
 \nonumber
 \int_t^{t+T(\Delta,\delta)} x^2_f(\tau) d\tau &\geqslant   \bar x^{2}_f T^\varepsilon-\psi(\Delta,\delta)\eta T^\varepsilon=\frac{\bar x^{2}_f}{2}T^\varepsilon.
 \end{align}
Thus, condition $(i)$ of item 2) of Lemma \ref{lem_PE} holds with $\mu(\Delta,\delta)=\frac{\bar x^{2}_f}{2}T^\varepsilon$. $\hfill\blacksquare$

{\bf Proof of Theorem \ref{thm_1}.} 
In the following, we denote the $(x_f,x_s)$-component of the solutions to (\ref{3D}) as $x$. We first prove item 1). Let $\Delta>\delta>0$, $\beta<\beta_c$ and constant input $\delta<\vert u\vert<\Delta$, according to item 1) of Proposition 1, there exists $\bar \varepsilon>0$ such that for any $\varepsilon\in(0,\bar\varepsilon]$ and any $\vert x(0)\vert \leqslant \Delta$, the corresponding solution $x$ to subsystem (\ref{3D1_1})-(\ref{3D1_2}) satisfies for $t\geqslant 0$,
\begin{equation}
\label{subsys_d}
\vert x(t)-x^*\vert\leqslant c_{1}(\Delta,\delta)e^{-c_{2}(\Delta,\delta)t}
\vert x(0)-x^*\vert,
\end{equation}
where $x^*:=(x^*_f,x^*_s)$ is the fixed point of subsystem (\ref{3D1_1})-(\ref{3D1_2}), and $c_{1}(\Delta,\delta), c_{2}(\Delta,\delta)>0$. 

We next consider for $t\geqslant 0$, 
\begin{equation}
\label{lya_auxi2}
W(t,\hat\beta):=p(t)(\hat\beta-\bar{\hat\beta})^2,
\end{equation}
where $\bar{\hat\beta}$ is the equilibrium point of (\ref{3D1_3}) associated with constant input $(x^*,u)$, and $p(t)\!\!:=\!\!\int_t^\infty \Phi^2(s,t)\; ds$ as in the proof of Proposition \ref{proposition2}. For any $\vert \hat\beta(0)\vert\leqslant \Delta$, by the similar computations in (\ref{phi_up})-(\ref{phi_down}) for any $t\geqslant 0$
\begin{eqnarray}
\label{bound_lyapunow_thm}
\underline M(\Delta,\delta) \vert \hat\beta-\bar{\hat\beta}\vert^2&\leqslant& W(t,\hat\beta)\leqslant \overline M(\Delta,\delta)\vert \hat\beta-\bar{\hat\beta}\vert^2,
\end{eqnarray}
where $\underline{M}(\Delta,\delta), \overline{M}(\Delta,\delta)>0$. Following similar lines as in (\ref{dP_fin}), we obtain $\dot p(t)-2kx^2_{f}p(t)=-1$, and
we derive that for any $t\geqslant 0$, where the time arguments are omitted,

\begin{eqnarray}
\nonumber
&&\frac{\partial W}{\partial t}+\frac{\partial W}{\partial \hat\beta} \hat f(\hat\beta, x_{f}, u-x_{s})\\\nonumber
&\leqslant& -\vert \hat \beta-\bar{\hat\beta}\vert^2+\nu(\Delta,\delta)\vert \hat\beta-\bar{\hat\beta}\vert\left(\vert x_{f}-x^*_{f}\vert+\vert x_{s}-x^*_{s}\vert\right),\\\label{dw2_2}
\end{eqnarray}
where $\nu(\Delta,\delta)>0$.
We deduce from the above inequality
\begin{eqnarray}
\nonumber
\vert \beta(t)-\bar{\hat\beta}\vert^2&\leqslant&c_3(\Delta,\delta) e^{-c_4(\Delta,\delta) t}\vert \hat \beta(0)-\bar{\hat\beta}\vert^2\\
&&+c_5(\Delta,\delta) e^{-c_4(\Delta,\delta)t}\vert x(0)- x^*\vert^2,
\end{eqnarray}
with $c_3(\Delta,\delta), c_4(\Delta,\delta), c_5(\Delta,\delta)>0$.
This property together with (\ref{subsys_d}) imply that the semiglobal exponential stability of system (\ref{3D}), which, in turn, implies the global asymptotic stability, according to \cite[Proposition 3.4]{angeli2002}. Thus, the desired result follows.

We next prove item 2). Let $\Delta>0$. For any $0<\beta-\beta_c<1$, according to item \rm{2)} of Proposition 1, there exist $\bar \varepsilon>0,\bar u\in(0,\Delta)$ such that for any $\varepsilon\in(0,\bar\varepsilon]$, for any $u\in(-\bar u,\bar u)$, there exists $m>0$ such that for any $\vert x(0)\vert_{P^\varepsilon} \leqslant m$ except the unique unstable fixed point of (\ref{3D1_1})-(\ref{3D1_2}), the corresponding solution $x$ to subsystem (15a)-(15b) satisfies for $t\geqslant 0$,
\begin{equation}
\label{subsys_d2}
\vert x(t)\vert_{P^\varepsilon}\leqslant c_6(m)e^{-c_7(m)t}
\vert x(0)\vert_{P^\varepsilon},
\end{equation}
where $P^\varepsilon $ is the periodic solution of subsystem (\ref{3D1_1})-(\ref{3D1_2}) and $c_6(m),c_7(m)>0$. By virtue of Lemma 2 and Proposition 2, subsystem (\ref{3D1_3}) is semiglobally incrementally ISS. Then according to \cite[Proposition 4.4] {angeli2002}, there exists initial condition $\bar{\hat\beta}^{lc}(0)$ such that $\bar{\hat\beta}^{lc}(t)$ is a periodic solution to (\ref{3D1_3}) with periodic input $(P^\varepsilon,u)$. We choose $W(t,\hat\beta):=p(t)(\hat\beta-\bar{\hat\beta}^{lc})^2$. Then using similar arguments as above we obtain the local exponential stability of system (\ref{3D}) with respect to a limit cycle.

Due to \cite[Proposition 4.5]{angeli2002}, for any initial condition $x(0)$ except the unique unstable fixed point $x^*$ of subsystem (\ref{3D1_1})-(\ref{3D1_2}), the corresponding solution to system (\ref{3D1_3}) holds $\lim\limits_{t\to\infty}\vert \hat\beta(t)\vert_{\bar{\hat\beta}^{lc}}=0$ as $\lim\limits_{t\to\infty}\vert x(t)\vert_{P^\varepsilon}=0$ . In addition, local exponential stability implies local stability, then the desired result follows.

{\bf Proof of Theorem \ref{thm_betah}.} We start by proving item \rm{1)}. By Theorem \ref{thm_1}, if $\beta<\beta_c$ then $(x_f(t),x_s(t))\to(x^*_f,x^*_s)$ and $\hat\beta(t)\to \bar{\hat\beta}$ as time goes to infinity. Note that, since $(x^*_f,x^*_s)$ is a steady-state of (\ref{3D1_1})-(\ref{3D1_2}), it must satisfies $(x^*_f,x^*_s)\in S^0$. Moreover, $u\neq 0$ implies $x^*_f\neq 0$. Because $-kx^*_f(-x^{*3}_f+\bar{\hat\beta} x^*_f+u-x^*_s)=0$ and $(x^*_f,x^*_s)\in S^0$, it follows that $\bar{\hat\beta}=\frac{x^{*3}_f-S^{-1}(x^*_f)+\beta x^*_f}{x^*_f}=\hat\beta^*(x^*_f,\beta)$.

\begin{figure}
\centering
\includegraphics[width=0.75\columnwidth]{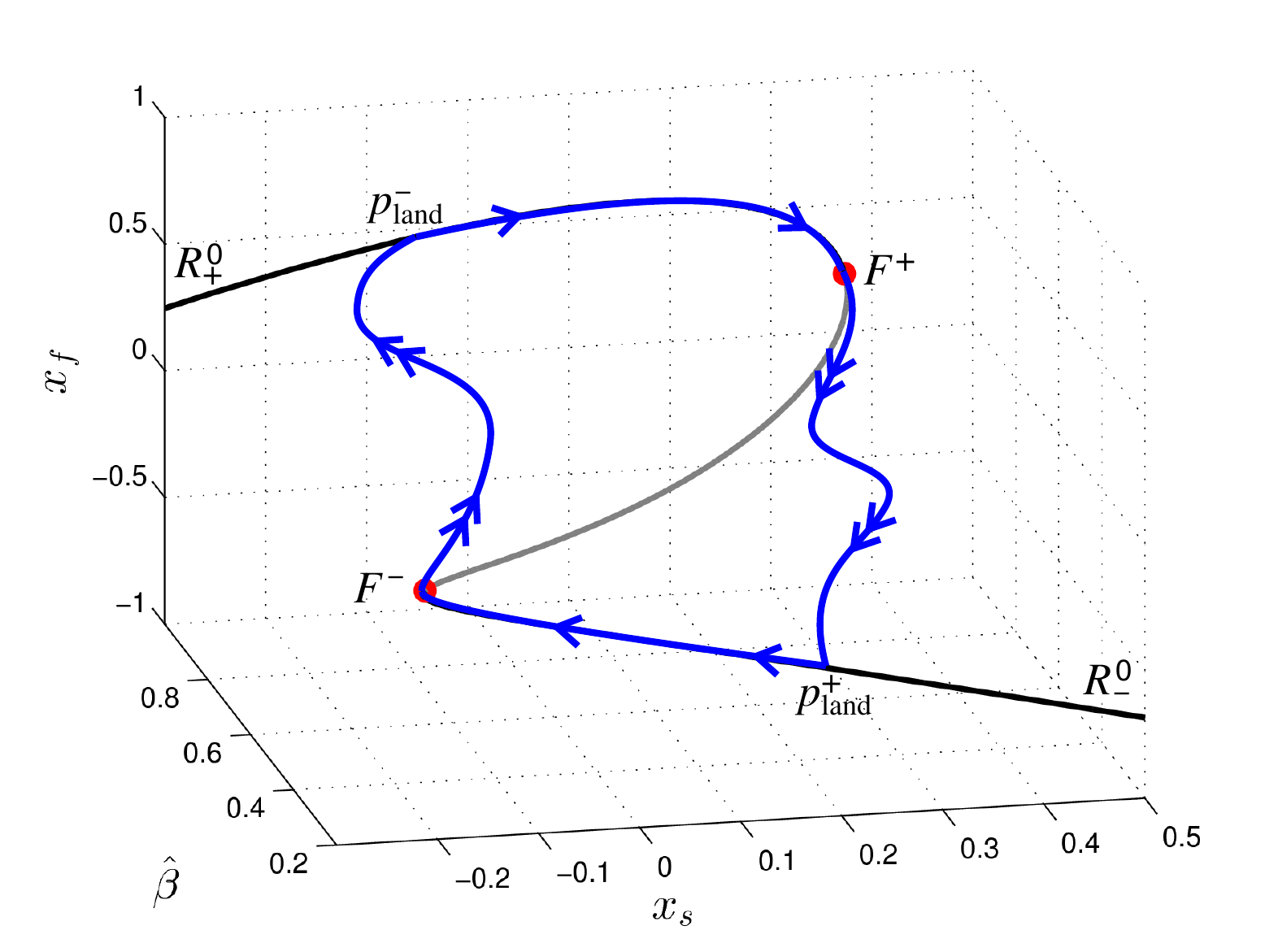}
\caption{}\label{FIG: 3D LC singular}
\end{figure}

We next prove item \rm{2)}. At the singular limit $\epsilon=0$, the one dimensional estimate critical manifold $R^0$ of system (\ref{3D}) is defined by the following two equations
\begin{subequations}\nonumber
\begin{align}[left = \empheqlbrace\,]
\nonumber
0&=-x_f+S(\beta x_f+u-x_s),\\\nonumber
0&=-x_f^3+\hat\beta x_f+u-x_s,
\end{align}
\end{subequations}
which can explicitly be solved as
\begin{IEEEeqnarray}{rCl}\label{EQ: obs crit man}
\IEEEyesnumber
R^0&=&\Bigg\{(x_f,x_s,\hat\beta):\\
&&x_s=\phi(x_f,u,\beta):=-S^{-1}(x_f)+\beta x_f+u,\nonumber\\
 &&\hat\beta= \frac{x_f^3-u+\phi(x_f,u,\beta)}{x_f}=\hat\beta^*(x_f,\beta)\Bigg\}.\nonumber
\end{IEEEeqnarray}
Note that the $(x_f,x_s)$-projection of $R^0$ is $S^0$.\\
The layer dynamics reads
\begin{subequations}\label{EQ: est layer dynamics}
\begin{align}[left = \!\!\!\!\empheqlbrace\,]
\dot x_f&=-x_f+S(\beta x_f+u-x_s),\IEEEyessubnumber\\
\dot x_s&=0,\IEEEyessubnumber\\
 \dot{\hat\beta}&=-kx_f(-x_f^3+\hat\beta x_f+u-x_s).\IEEEyessubnumber
\end{align}
\end{subequations}
To study the normal hyperbolicity of $R^0$, we compute the Jacobian matrix of (\ref{EQ: est layer dynamics}) on $R^0$ as 
\begin{eqnarray}
\nonumber
\!\!\!\!J&\!\!=\!\!&\begin{pmatrix}
\!\!\!\!\!\!-1+\beta S^{\prime}(\beta x_f+u-x_s)  &  \!\!\!\!-S^{\prime}(\beta x_f+u-x_s) &0  \\
0& 0&0\\
4k x_f^3-2k\hat\beta x_f -k(u-x_s)  & kx_f& -k x_f^2
\end{pmatrix}
\end{eqnarray}
for $x_s=\phi(x_f,u,\beta)$, $\hat\beta=\hat\beta^*(x_f,\beta)$.
The two nonzero eigenvalues are $\lambda_1=-kx^2_f,\;\lambda_2=-1+\beta S^\prime(S^{-1}(x_f))$.
For $0<\beta-\beta_c<1$, it holds that $\beta>\beta_c=\frac{1}{S^\prime(0)}$, which is equivalent to $\frac{1}{\beta}<S^\prime(0)$ and therefore, invoking item d) of Assumption \ref{sigmo_assum}, the equation $-1+\beta S^\prime(S^{-1}(x_f))=0$ has exactly one positive and negative roots, $x_{fold}^{+}$ and $x_{fold}^{-}$, respectively. Moreover,
\begin{equation*}
\lambda_2\left\{\begin{array}{rl}
<0, & {\rm if}\ x_f<x_{fold}^{-}\ {\rm or}\ x_f>x_{fold}^{+},\\
>0, & {\rm if}\ x_{fold}^{-}<x_f<x_{fold}^{+}.
\end{array}\right.
\end{equation*}
Since $\lambda_1<0$ for $x_f\neq 0$, it follows that two branches of $R^0$ given by
\begin{IEEEeqnarray*}{rCl}
R^0_+&:=&R^0\cap\{x_f>x_{fold}^{+}\},\\
R^0_-&:=&R^0\cap\{x_f<x_{fold}^{-}\},
\end{IEEEeqnarray*}
are locally exponentially attractive branches of the estimate critical manifold $R^0$. Let
$$\hat\beta_{fold}^+=\hat\beta^*(x_{fold}^+,\beta),\; \hat\beta_{fold}^-=\hat\beta^*(x_{fold}^-,\beta),$$
$$ x_{s,fold}^+=\phi(x_{fold}^+,u,\beta) 
,\; x_{s,fold}^-=\phi(x_{fold}^-,u,\beta).$$
The end points of $R^0_+$ and $R^0_-$ are given by
$$F^+:=(x_{fold}^+,\hat\beta_{fold}^+,x_{s,fold}^+),\; F^-:=(x_{fold}^-,\hat\beta_{fold}^-,x_{s,fold}^-) $$
respectively, are fold points of the estimate critical manifold $R^0$. See Figure~\ref{FIG: 3D LC singular}.

The slow flow on $R^0$ is defined by the reduced dynamics
\begin{subequations}\label{EQ: est reduced dynamics}
\begin{align}[left = \empheqlbrace\,]
0&=-x_f+S(\beta x_f+u-x_s),\IEEEyessubnumber\\
x_s^\prime&=x_f-x_s,\IEEEyessubnumber\\
0&=-kx_f(-x_f^3+\hat\beta x_f+u-x_s).\IEEEyessubnumber
\end{align}
\end{subequations}

The flow on $R^0_+$, as well as on $R^0_-$, with respect to the reduced dynamics of (\ref{EQ: est reduced dynamics}) is given as follows. It holds that $x_s^\prime <0 (x_f<x_s<0)$ on $R^0_-$ and, similarly, $x_s^\prime>0$ on $R^0_+$. Noticing that the $x_s$-projection of $R^0_+$ (resp. $R^0_-$) is the semiline $(-\infty, x^+_{s,fold})$ (resp. $(x^-_{s,fold}, \infty)$), it follows that all trajectories on $R^0_+$ eventually reach the fold point $F^+$. Conversely, all trajectories on $R^0_-$ eventually reach $F^-$ (see Fig.~\ref{FIG: 3D LC singular}).

At the folds, we connect the slow flow with the (two-dimensional) fast flow. We claim that the fast flow brings the trajectory on the opposite branch of the critical manifold. Indeed, by monotonicity of (\ref{EQ: est layer dynamics}a), $x_f$ converges to equilibrium. By the cascade structure of (\ref{EQ: est layer dynamics}a)-(\ref{EQ: est layer dynamics}c) and incremental ISS property of (\ref{EQ: est layer dynamics}c), $\hat\beta$ also converges to equilibrium. Because the only equilibria of (\ref{EQ: est layer dynamics}a)-(\ref{EQ: est layer dynamics}c) are on the estimate critical manifold, the result follows.

We have therefore constructed a candidate singular periodic orbit $ Q^0$. The persistence of this singular orbit for $\varepsilon>0$ follows as in the proof of Proposition \ref{LEM: stable osci dichotomy} and is omitted. Let $Q^\varepsilon$ be the resulting periodic orbit of (\ref{3D}). 
By Theorem \ref{thm_1}, almost all the trajectories of system (\ref{3D}) asymptotically converge to $Q^\varepsilon$. Following again the same arguments as the proof of Proposition \ref{LEM: stable osci dichotomy}, it holds that trajectories along $Q^\varepsilon$ spend only an $\mathcal O(\varepsilon)$-fraction of the limit cycle period outside an $\mathcal O(\varepsilon)$ neighborhood of the estimate critical manifold. Because $\hat\beta(t)=\hat\beta^*(x^{lc}_f(t+\theta),\beta)$ on the estimate critical manifold, the result follows.$\hfill\blacksquare$

\bibliographystyle{unsrt}
\bibliography{tyartr}

\end{document}